\documentclass[journal]{ieeetran}
\usepackage{cite}
\usepackage{graphicx}
\usepackage{epstopdf}
\usepackage{amsmath}
\usepackage{cases}
\usepackage[caption=false,font=footnotesize]{subfig}
\usepackage{fixltx2e}
\usepackage{amssymb}
\usepackage{mathtools}
\usepackage{bm}
\usepackage{multirow}
\usepackage{color}
\usepackage{algorithmic}
\usepackage[linesnumbered, ruled, vlined]{algorithm2e}
\usepackage{amsthm}
\usepackage{mathrsfs}
\usepackage{textcomp,mathcomp}
\usepackage{threeparttable}

\newcommand{\tabincell}[2]{\renewcommand\arraystretch{0.9}\begin{tabular}{@{}#1@{}}#2\end{tabular}}

\begin{document}

\title{Optimal Investment Portfolio of Thyristor- and IGBT-based Electrolysis Rectifiers in\\Utility-scale Renewable P2H Systems}
\author{
Yangjun~Zeng,~\IEEEmembership{Student Member,~IEEE},
Yiwei~Qiu,~\IEEEmembership{Member,~IEEE},
Liuchao~Xu,
Chenjia~Gu,~\IEEEmembership{Member,~IEEE},
Yi~Zhou,~\IEEEmembership{Member,~IEEE},
Jiarong~Li,~\IEEEmembership{Member,~IEEE},
Shi~Chen,~\IEEEmembership{Member,~IEEE}, and
Buxiang~Zhou,~\IEEEmembership{Member,~IEEE}

\thanks{Financial support was obtained from the National Key R\&D Program of China (2021YFB4000503), the National Natural Science Foundation of China (52377116 and 52307126), and the Natural Science Foundation of Sichuan Province (2024NSFSC0870). \emph{(Corresponding author: Yiwei Qiu)}}
\thanks{Y. Zeng, Y. Qiu, L. Xu, C. Gu, Y. Zhou, S. Chen, and B. Zhou are with the College of Electrical Engineering, Sichuan University, Chengdu 610065, China.
(ywqiu@scu.edu.cn)
}
\thanks{J. Li is with the Harvard John A. Paulson School of Engineering and Applied Sciences, Harvard University, Cambridge 02138, USA.

}%
}
\maketitle
\begin{abstract}
  Renewable power-to-hydrogen (ReP2H) systems require rectifiers to supply power to electrolyzers (ELZs). Two main types of rectifiers, insulated-gate bipolar transistor rectifiers (IGBT-Rs) and thyristor rectifiers (TRs), offer distinct tradeoffs. IGBT-Rs provide flexible reactive power control but are costly, whereas TRs are more affordable with lower power loss but consume a large amount of uncontrollable reactive power. A mixed configuration of rectifiers in utility-scale ReP2H systems could achieve a decent tradeoff and increase overall profitability. To explore this potential, this paper proposes an optimal investment portfolio model. First, we model and compare the active and reactive power characteristics of ELZs powered by TRs and IGBT-Rs. Second, we consider the investment of ELZs, rectifiers, and var resources and coordinate the operation of renewables, energy storage, var resources, and the on-off switching and load allocation of multiple ELZs. Subsequently, a two-stage stochastic programming (SP) model based on weighted information gap decision theory (W-IGDT) is developed to address the uncertainties of the renewable power and hydrogen price, and we apply the progressive hedging (PH) algorithm to accelerate its solution. Case studies demonstrate that optimal rectifier configurations increase revenue by at most {\color{black}13.78\% compared with configurations using only TRs or IGBT-Rs, existing project setups, or intuitive designs}. Under the optimal portfolio, reactive power compensation investment is nearly eliminated, with a preferred TR-to-IGBT-R ratio of 3:1.
\end{abstract}


\section*{Nomenclature}
{\subsection{Abbreviations}
\begin{IEEEdescription}[\IEEEusemathlabelsep\IEEEsetlabelwidth{superscript}]
\addcontentsline{toc}{section}{Nomenclature}
\item[BES] Battery energy storage
\item[ELZ] Electrolyzer
\item[IGBT-R] Insulated-gate bipolar transistor-based rectifier
\item[OLTC] On-load tap changer
\item[PV] Photovoltaic
\item[ReP2H] Renewable power-to-hydrogen
\item[SP] Stochastic programming
\item[TR] Thyristor rectifier
\item[W-IGDT] Weighted information gap decision theory
\item[WT] Wind turbine
\end{IEEEdescription}}

\subsection{Indices and Sets}
\begin{IEEEdescription}[\IEEEusemathlabelsep\IEEEsetlabelwidth{superscript}]
\addcontentsline{toc}{section}{Nomenclature}

\item[$n,t$] Index for ELZs and time periods
\item[$i, j,j'$] Index for buses
\item[$ij$]   Index for branches
\item[$\pi(j)$, $\sigma(j)$] Set of parents and children of bus $j$

\end{IEEEdescription}

\subsection{Variables}

\subsubsection{Investment Variables}
\begin{IEEEdescription}[\IEEEusemathlabelsep\IEEEsetlabelwidth{superscript}]
\addcontentsline{toc}{section}{Nomenclature}
\item[$\sigma^{\text{ELZ}}$] Number of ELZs
\item[$\sigma^{\text{TR}}, \sigma^{\text{IGBT}}$] Number of TR and IGBT-R
\item[$W^{\text{C}}$] Capacity of var compensation
\end{IEEEdescription}

\subsubsection{ELZ-related Operational Variables}
\begin{IEEEdescription}[\IEEEusemathlabelsep\IEEEsetlabelwidth{superscript}]
\addcontentsline{toc}{section}{Nomenclature}

\item[$P^{\text{Stack}}, U^{\text{Stack}}$] DC power and voltage of the stack
\item[$I, Y^{\text{H}_2}$] Electrolytic current and hydrogen flow
\item[$T$] Stack temperature of the ELZ
\item[$b^{\text{On}}, b^{\text{By}}, b^{\text{Idle}}$] Production, standby, and idle states
\item[$b^{\text{SU}}, b^{\text{SD}}$] Startup/shutdown actions of ELZ
\item[$P^{\text{ELZ}}, P^{\text{Rec}}$] Active power of the ELZ/rectifier
\item[$P^{\text{Loss}}$] The active loss of the TR/IGBT-R of the ELZ
\item[$Q^{\text{TR}}, Q^{\text{IGBT}}$] Reactive power of the TR and IGBT-R
\item[$P^{\text{BoP}}$] Power of the balance of plant of the ELZ
\item[$P^{\text{Heat}}, P^{\text{Diss}}$] Electrolytic heat and dissipation of the ELZ
\item[$P^{\text{Cool}}$] Cooling heat flow of the ELZ
\end{IEEEdescription}

\subsubsection{Network-side Operational Variables}
\begin{IEEEdescription}[\IEEEusemathlabelsep\IEEEsetlabelwidth{superscript}]
\addcontentsline{toc}{section}{Nomenclature}
\item[$P_{ij,t}, Q_{ij,t}$] Active/reactive power flows on branch $ij$
\item[$I_{ij,t}, \ell_{ij,t}$] Current on branch $ij$ and its square
\item[$U_{j,t}, \upsilon_{j,t}$] Voltage amplitude of bus $j$ and its square
\item[$p_{j,t}, q_{j,t}$] Active/reactive power injections at bus $j$ \vspace{1pt}
\item[$p_{j,t}^{\text{L}}, q_{j,t}^{\text{L}}$] Active/reactive loads of all ELZs at bus $j$ \vspace{1pt}
\item[$P_{j,t}^{\text{WT}}, Q_{j,t}^{\text{WT}}$] Active/reactive power of the WT at bus $j$ \vspace{1pt}
\item[$P_{j,t}^{\text{PV}}, Q_{j,t}^{\text{PV}}$] Active/reactive power of the PV at bus $j$ \vspace{1pt}
\item[$E_{j,t}$]  Energy stored level of the BES at bus $j$ \vspace{1pt}
\item[$P_{j,t}^{\text{BES,C/D}}$] BES charging/discharging power at bus $j$ \vspace{1pt}
\item[$b_{j,t}^{\text{BES,C/D}}$] BES Charging/discharging states at bus $j$ \vspace{1pt}
\item[$Q_{j,t}^{\text{BES/C}}$] Reactive power of the BES/var comp. at bus $j$
\item[$k_{ij,t},\delta_{ij,k,t}$] Turn ratio of the OLTC and its binary decision on the $k$th level of the transformer $ij$
\end{IEEEdescription}

\subsection{Parameters}

\begin{IEEEdescription}[\IEEEusemathlabelsep\IEEEsetlabelwidth{superscript}]
\addcontentsline{toc}{section}{Nomenclature}

\item[${N}_s$] The product of the number of typical days and the number of hydrogen price scenarios
\item[${N_t}, \Delta t$] Scheduling horizon and step length
\item[$c^{\text{ELZ}}$] Investment cost of ELZ
\item[$c^{\text{TR}},c^{\text{IGBT}}$] Investment costs of TR and IGBT-R
\item[$c^{\text{C}}$] Investment cost of var compensation
\item[$c^{\text{H}_2}$] Hydrogen price
\item[$c^{\text{SU}}, c^{\text{SD}}$] Startup/shutdown costs of the ELZ
\item[$\overline{\sigma}^{\text{ELZ}}$] Upper limit of the number of ELZs
\item[$\overline I, \underline I$] Electrolytic current limits of the ELZ
\item[$\overline T, \underline T$] Stack temperature limits of the ELZ
\item[$\eta^{\text{Cool}}$] Cooling efficiency of the ELZ
\item[$P^{\text{By}}$] Standby power consumption of the ELZ
\item[$C^{\text{ELZ}}, R^{\text{Diss}}$] Heat capacity and dissipation resistance of ELZ
\item[$T^{\text{Am}}$] Ambient temperature
\item[$c^{\text{Cool}}, T^{\text{Cool}}$] Cooling factor and coolant temperature of ELZ
\item[$r_{ij}, x_{ij}$] Resistance and reactance of branch $ij$
\item[$\overline{U}_j, \underline{U}_j$] Voltage magnitude limits at bus $j$
\item[$\overline{I}_{ij}$] Current capacity limit of branch $ij$
\item[$S_{j}^\text{WT/PV/BES}$] WT/PV/BES installation capacities at bus $j$\vspace{1pt}
\item[$S^\text{sc}$] Short-circuit capacity at hydrogen plant node
\item[$S^\text{TR/IGBT}$] Rated capacities of TR/IGBT-R
\item[$\theta$] Minimum power factor angle of PV plants
\item[$\overline{P}^{\text{BES,C/D}}$] BES charging/discharging power limits
\item[$\eta^{\text{BES,C/D}}$] BES charging/discharging efficiencies
\item[$\zeta^{\text{BES}}$] BES self-discharge ratio
\item[$\overline{E}, \underline{E}$] Energy stored level limits of the BES
\item[$\overline{k}_{ij}^{\text{All}}$] OLTC switching time limit of transformer $ij$
\item[$\overline{k}_{ij},\underline{k}_{ij}$] OLTC ratio limits of transformer $ij$
\end{IEEEdescription}

\section{Introduction}
\label{sec:intro}

\subsection{Background and Motivation}

\IEEEPARstart{R}{enewable} power-to-hydrogen (ReP2H) facilitates wind and solar integration \cite{van2020hydrogen} and has significant potential for applications in energy storage, transportation \cite{van2020hydrogen}, and the chemical industry \cite{guo2023deploying}. With the continuous growth of renewable energy (RE) installations, the ReP2H industry has grown rapidly under policy incentives and market drivers \cite{odenweller2022probabilistic,2022}, leading to increased investment in electrolyzers (ELZs) in large-scale hydrogen plants to replace fossil-based hydrogen production.

Industrial hydrogen plants typically consist of multiple ELZs \cite{qiu2023extend}, each integrating a rectifier, an electrolytic stack, and the balance of plant (BoP) \cite{zeng2024scheduling}. Since electrolytic stacks operate on DC power, rectifiers are required to convert AC to DC. Currently, two types of rectifiers are commercially available {\color{black}and widely used in engineering}: thyristor rectifiers (TRs) and  {\color{black}voltage-source} insulated-gate bipolar transistor-based rectifiers (IGBT-Rs) \cite{koponen2021comparison, songyuan2023, daan2024}. They offer distinct tradeoffs, as summarized in Table \ref{tab:rectifier}.
{\color{black}Besides these, diode-based rectifiers (DRs), current-source rectifiers (CSRs) \cite{wu2024coordination}, and hybrid thyristor-IGBT rectifiers (HRs) \cite{meng2022novel} have also been studied. However, due to the limited controllability of DRs, the immaturity of the CSRs \cite{wu2024coordination}, and the engineering complexity of HRs, these rectifiers are not commercially available on large scales. Hence, only SCRs and IGBT-Rs are adopted in utility-scale ReP2H projects, and this work focuses on these two types.}

IGBT-Rs use PWM modulation and voltage vector control, enabling independent active and reactive power control. Thus, IGBT-Rs can maintain a power factor (PF) of 1 \cite{koponen2021comparison} or provide reactive power support \cite{de2023hydrogen}, or even enable low voltage ride-through \cite{zhang2022research}. Moreover, their fast response makes IGBT-Rs well suited for weak or off-grid systems \cite{tavakoli2023gridforming}. However, their structural complexity leads to higher costs and power losses.

TRs are more mature and widely applied in low-voltage high-current applications such as water and chlor-alkali electrolysis. They employ multipulse (usually 12--48) rectification to reduce harmonics. They are more cost-effective and efficient than IGBT-Rs.
However, owing to their phase-controlled nature, TRs consume a large amount of reactive power, and the active and reactive power components are coupled \cite{ruuskanen2020power}. Thus, additional var resources, such as compensators and on-load tap changers (OLTCs), are often required to ensure voltage security \cite{zeng2024scheduling,ruuskanen2020power}.

In practice, some large-scale projects \cite{songyuan2023,daan2024} use mixed rectifier (MR) configurations to leverage the complementary characteristics of TRs and IGBT-Rs. The \emph{China Energy Engineering Songyuan Hydrogen Industry Park} \cite{songyuan2023} in Songyuan, China, employs 32 TRs and 32 IGBT-Rs to power 64 alkaline ELZs, and the \emph{Wind-Solar Hydrogen-Ammonia Integration Project} \cite{daan2024} in Da'an, China, uses 12 TRs and 24 IGBT-Rs for 36 ELZs. However, to the best of the author's knowledge, such configurations have not been quantitatively analyzed.
{\color{black}Overlooking rectifier performance in the planning stage may hinder coordination between the hydrogen plant and the electrical network, as the impact of reactive power loads and resources on ReP2H system economics becomes hard to quantify.}
Suboptimal configurations may lead to excessive or insufficient reactive power compensation{ \color{black}and reduced hydrogen yield}, increasing investment or operating costs. {\color{black}Illustrative examples are presented in Section \ref{sec:comparision}.}

{\color{black}
	Specifically, the optimal investment of rectifiers in a ReP2H system is a comprehensive optimization problem rather than a simple pairing of ELZs with rectifiers. It needs multiple interrelated considerations:
	\begin{enumerate}
		\item  The coupling of active and reactive power in TR-powered ELZs, which presents tradeoffs between improving P2H efficiency and managing reactive power demand \cite{zeng2024scheduling,li2024two}.
		\item  Coupling between active and reactive power flows in the electrical network of the ReP2H system.
		\item  Coupling between the hydrogen production process and the network power flow.
	\end{enumerate}
}
However, existing ReP2H planning research has focused mainly on the sizing of energy sources and ELZs while considering only active power and neglecting rectifier selection, {\color{black}the electrical network}, and var resources. {\color{black}Coordinating investment planning and production scheduling with respect to both active and reactive power} is vital for improving capital efficiency, which, unfortunately, has not been investigated.

\begin{table}[tb]\scriptsize
	\renewcommand{\arraystretch}{1.7}
	\caption{Comparison of TR and IGBT-R Powering a 5-MW Electrolyzer}\vspace{-0pt}
	\label{tab:rectifier}
	\centering
	\begin{threeparttable}
		\begin{tabular}{c@{\hspace{5pt}}c@{\hspace{5pt}}c@{\hspace{5pt}}c@{\hspace{5pt}}c}
			\hline\hline
			Rectifier             & Reactive power characteristics    & \tabincell{c}{Investment\\cost (CNY)}  & \tabincell{c}{Efficiency}         \\
			\hline
			TR\footnote           & \tabincell{c}{Coupled with active load \cite{ruuskanen2020power}; see Fig. \ref{fig:region}(a)}    & $\approx$1 million   & \tabincell{c}{$\approx$98\% \cite{gao2024advanced}}        \\
			\hline
			IGBT-R                & \tabincell{c}{PF=1 \cite{koponen2021comparison}; or to provide and absorb\\ reactive power decoupled from\\ active load \cite{de2023hydrogen,zhang2022research}; see Fig. \ref{fig:region}(b)}    & 2 to 3 million    & \tabincell{c}{$\approx$97\% \cite{gao2024advanced}}           \\
			\hline\hline
		\end{tabular}
		\begin{tablenotes}
			\footnotesize
			\item[1] Experiments \cite{gao2024advanced,intellipower2022} have shown that a 24-pulse TR can control total harmonic distortion (THD) of current to 3--5\% (the IGBT-R is even better), which can meet the power grid regulation with filters. Therefore, this study does not consider the harmonics of either type and leaves this topic to future work.
		\end{tablenotes}
	\end{threeparttable}
	\vspace{-0pt}
\end{table}

\subsection{Literature Review and Research Gap}

ReP2H planning research typically addresses RE, battery energy storage (BES), ELZs, and hydrogen storage (HS). For example, Li \textit{et al.} \cite{li2019optimal} proposed an investment model for ELZs and seasonal HS, considering power and hydrogen delivery, to address the spatiotemporal imbalance between renewable resources and hydrogen demand. Pan \textit{et al.} \cite{pan2020optimal} developed a robust capacity sizing model for wind turbines (WTs), photovoltaics (PVs), ELZs, and HSs in integrated energy systems. Zhu \textit{et al.} \cite{zhu2024full} proposed an energy management framework and studied BES capacity sizing in off-grid ReP2H systems.
Li \textit{et al.} \cite{li2023exploration} explored the configuration of the number and size of ELZs and developed operation rules to enhance wind power integration. To assess the cost competitiveness of green hydrogen, Ib{\'a}{\~n}ez-Rioja \textit{et al.} \cite{ibanez2023off} optimized planning and operation using 30 years of real-world wind and solar data. Similarly, Zheng \textit{et al.} \cite{zheng2023model} analyzed the levelized cost of hydrogen (LCOH) for off-grid wind P2H systems and reported that optimal sizing reduces LCOH. However, the aforementioned studies \cite{li2019optimal,pan2020optimal,zhu2024full,li2023exploration,ibanez2023off,zheng2023model} overlooked rectifier configuration and reactive power balance in ReP2H planning. {\color{black}This omission could potentially cause ineffective var resource sizing from disregarded rectifier reactive characteristics and reduce hydrogen productivity due to unconsidered rectifier losses. Moreover,} given that rectifiers account for 20--30\% of ELZ investment \cite{escn2024}, their optimal selection is crucial for economic feasibility.

{\color{black}The effectiveness of planning is closely associated with operational performance. Regarding P2H scheduling, Wang \textit{et al.} \cite{wang2024optimization} proposed a two-phase control strategy that considers the minimum power constraint of ELZs, determining the on-off switching and power to maximize RE utilization while balancing their lifespan.
Liang \textit{et al.} \cite{liang2024large} designed an ELZ scheduling strategy using the Pelican Optimization Algorithm to improve P2H efficiency. 
Firdous \textit{et al.} \cite{firdous2025short} developed a modular P2H management model in integrated power-to-ammonia systems to enhance operation flexibility. Khajeh \textit{et al.} \cite{Khajeh2024Optimized} proposed a wind-hydrogen system scheduling scheme to provide frequency support and congestion management for the transmission system, thereby enhancing its profitability. Han \textit{et al.} \cite{Han2025Robust} proposed a robust optimization framework for an integrated electricity-heat-hydrogen system that accounts for ELZ-network interactions. 
Even with their contributions, these studies did not account for rectifier reactive power or its coupling with the electrical network, making it difficult to quantify the impact of var resources on ReP2H system economics.
Based on rectifier performance, the active-reactive power coordination between hydrogen production and the network is particularly important.}

{\color{black}Current research on electrolysis rectifiers} has focused mainly on performance analysis \cite{ruuskanen2020power,koponen2021comparison} or short-term scheduling \cite{zeng2024scheduling,li2024two}.
For example, Ruuskanen \textit{et al.} \cite{ruuskanen2020power} investigated the power quality and reactive power characteristics of ELZs powered by TRs. 
Koponen \textit{et al.} \cite{koponen2021comparison} compared the tradeoffs of the TR and IGBT-R in terms of performance and investment in industrial applications. In our prior work \cite{li2024two,zeng2024scheduling}, Li \textit{et al.} \cite{li2024two} integrated rectifier power characteristics into ReP2H scheduling, revealing conflicts between productivity and the grid-side PF in grid-connected hydrogen plants, whereas Zeng \textit{et al.} \cite{zeng2024scheduling} addressed the tradeoffs between P2H energy conversion efficiency and network losses caused by the reactive power of a TR and proposed a coordinated power management approach for large hydrogen plants that integrate multiple ELZs.

Despite these advancements, the optimal rectifier configuration in ReP2H planning remains unaddressed. Current designs are often based on intuition \cite{songyuan2023,daan2024} rather than quantitative analysis. The challenge of leveraging the complementarities between TRs and IGBT-Rs to achieve an optimal tradeoff remains unexplored.

\subsection{Contributions of This Work}

To fill the aforementioned gap and explore MR configurations for improving capital efficiency in ReP2H systems, we propose an optimal portfolio model based on stochastic programming (SP). The key contributions are as follows:

\begin{enumerate}
	\item The active and reactive characteristics of TR- and IGBT-R-powered ELZs are modeled and compared, {\color{black}demonstrating their complementarity.
		Further analysis shows that an MR configuration enables more effective coordination between hydrogen production and network-side performance in ReP2H systems, as compared to strategies focused solely on hydrogen plant operation.
	} 
	\item  An optimal investment portfolio model for MR and var resources in ReP2H systems is proposed, {\color{black}in which the operation stage adopts a coordinated active-reactive power scheduling strategy for hydrogen production and electrical network}. The model handles uncertainties via weighted information gap decision theory (W-IGDT). A two-stage decomposition algorithm, progressive hedging (PH) \cite{rockafellar1991scenarios}, is applied to improve computational efficiency.
	\item Case studies on realistic systems show that optimizing MR configuration can increase revenue by at most {\color{black}13.78\%} compared with that of uniform configurations or those in existing projects and {\color{black}intuitive designs}. The results suggest that a 3:1 ratio of TRs to IGBT-Rs offers better benefits {\color{black}and is expected to provide a reference for industrial practices}.
\end{enumerate}

The rest of this paper is organized as follows. Section \ref{sec:model} presents the load characteristics of ELZs powered by TRs and IGBT-Rs. Section \ref{sec:SOD} establishes the investment portfolio model and corresponding solution approach. In Section \ref{sec:cases}, simulations verify the proposed model. Finally, the conclusions drawn from this work are summarized in Section \ref{sec:conclusion}.

\section{System Structure and Comparison Between ELZs Powered by TR and IGBT-R}
\label{sec:model}


\subsection{General Configuration of ReP2H Systems}
\label{sec:ReP2H}

\begin{figure}[t]
	\centering
	\includegraphics[width=3.4in]{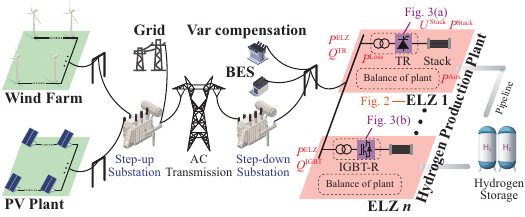}\vspace{-0pt}
	\caption{\color{black}The general structure of ReP2H systems.}
	\label{fig:system}\vspace{-0pt}
\end{figure}

Based on existing projects \cite{songyuan2023,daan2024}, the general structure of a ReP2H system is illustrated in Fig. \ref{fig:system}. The system integrates RE sources (wind and solar), transmission lines, a hydrogen plant, and sometimes a BES. The hydrogen plant often comprises multiple ELZs.
{\color{black} Due to policy restrictions \cite{neimenggu2024} and PF constraints \cite{li2024two}, which aim to ensure local renewable utilization, facilitate green hydrogen certification, and avoid occupying the transmission capacity, both active and reactive power exchange with the external grid is limited. Therefore, such exchanges are not considered in this study.}

The ELZs, as DC loads, may be powered by TRs or IGBT-Rs. They exhibit distinctive active and reactive load characteristics, as detailed in Section \ref{sec:reactpower}. This paper assumes that TRs are commonly used 24-pulse rectifiers \cite{gao2024advanced}.

\subsection{Hydrogen Production, Heat, and State-switching Characteristics of Electrolyzers}
\label{sec:same}
\begin{figure}[t]
	\centering
	\includegraphics[width=3.5in]{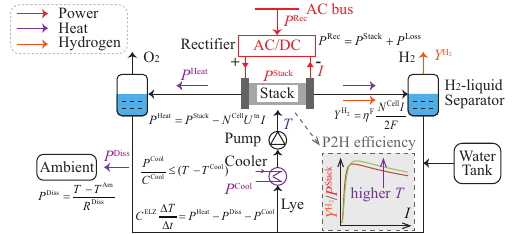}\vspace{-0pt}
	\caption{Schematic diagram of an alkaline electrolyzer.}
	\label{fig:ael}\vspace{-0pt}
\end{figure}

A schematic diagram of an alkaline ELZ is shown in Fig. \ref{fig:ael}. TR- and IGBT-R-powered ELZs share common characteristics, including hydrogen production, temperature dynamics, and state transitions. We summarize their physical models and operational constraints in Table \ref{sec:model}.

\subsubsection{Hydrogen production and heat submodels}

The active power of an ELZ consists of stack power, BoP consumption, and rectifier losses, modeled as (\ref{eq:Pelz})--(\ref{eq:Paux}). The electrolytic voltage-current relationship follows a nonlinear U-I curve \cite{ulleberg2003modeling}, given by (\ref{eq:Ustack}), with safety and operational limits of (\ref{eq:I}).

The power consumed by the electrolytic stack partially converts into hydrogen, while the remainder dissipates as heat. The hydrogen production flow $Y^{\text{H}_2}$ follows (\ref{eq:YH2}), where the Faradaic efficiency $\eta^{\text{F}}$ is modeled by (\ref{eq:etaF1}). The P2H energy conversion efficiency $Y^{\text{H}_2}/P^{\text{Stack}}$ is significantly affected by the temperature, as shown in Fig. \ref{fig:ael}.

The temperature control relies on the coordination between the stack and BoP. Proper thermal management of the ELZ can improve P2H efficiency and reduce auxiliary power use. We use a first-order model (\ref{eq:thermal}) to describe the temperature dynamics \cite{ulleberg2003modeling}, summarized as (\ref{eq:Pgen})--(\ref{eq:T}).
\begin{table}[tb]\scriptsize
	\renewcommand{\arraystretch}{1.4}
	\caption{Production, Heat, and State Transition Model of Electrolyzers}\vspace{-6pt}
	\label{tab:model}
	\centering
	\begin{tabular}{c@{\hspace{2pt}}c@{\hspace{5pt}}c}
		\hline\hline
		Submodel         & Physical model and operational constraints                  \\
		\hline
		\multirow{9}{*}{\tabincell{c}{Hydrogen\\production}}
		& \begin{minipage}{0.4\textwidth}
			\begin{equation}\text{ELZ power: }P^{\text{ELZ}}=P^\text{Rec} + P^\text{BoP} \label{eq:Pelz}
		\end{equation}\end{minipage}\\
		& \begin{minipage}{0.4\textwidth}\begin{equation}\text{Rectifier power: }P^{\text{Rec}}=P^\text{Stack} + P^\text{Loss},~P^{\text{Rec}}\eta^{\text{Rec}}=P^\text{Stack} \label{eq:Prec}
		\end{equation}\end{minipage}\\
		& \begin{minipage}{0.4\textwidth}\begin{equation}\text{Stack power: }P^{\text{Stack}}=U^{\text{Stack}}I  \label{eq:Pstack}
		\end{equation}\end{minipage}\\
		&\begin{minipage}{0.4\textwidth}\begin{equation}\text{Auxiliary power: }P^{\text{BoP}}=P^{\text{Cool}}/\eta^{\text{Cool}}+b^{\text{By}}P^{\text{By}}\label{eq:Paux}
		\end{equation}\end{minipage}\\
		& \begin{minipage}{0.4\textwidth} \vspace{-4pt} \begin{equation}\begin{split} \text{U-I characteristics: }
					U^{\text{Stack}}=N^{\text{Cell}} \big[ U^{\text{rev}}+({r_{1}+r_{2}T}){I}/{A} \\
					+ s_1\log\big({({t_{1}+{t_{2}}/{T}+{t_{3}}/{T^{2}}})}{I}/{A} + 1\big) \big]
				\end{split}\label{eq:Ustack}\end{equation}\end{minipage}\\
		&\begin{minipage}{0.4\textwidth}\begin{equation}\text{Current limits: }b^{\text{On}}\underline I\leq I\leq b^{\text{On}}\overline I  \label{eq:I}
		\end{equation}\end{minipage}\\
		& \begin{minipage}{0.4\textwidth}\begin{equation}\text{Hydrogen flow: }Y^{\text{H}_2}=\eta^{\text{F}}{N^{\text{Cell}}I}/(2F) \label{eq:YH2}
		\end{equation}\end{minipage}\\
		& \begin{minipage}{0.4\textwidth}\begin{equation}\text{Faradaic efficiency: }\eta^{\text{F}}={(I/A)^{2}}/[{{f_{1}+(I/A)^{2}}}]\times f_{2} \label{eq:etaF1}
		\end{equation}\end{minipage}\\
		\vspace{-9pt}
		\\
		\hline
		\multirow{5}{*}{Heat}
		& \begin{minipage}{0.4\textwidth}
			\vspace{-5pt} \begin{equation}\begin{split}\text{Thermal dynamics: }C_n^{\text{ELZ}}({T_{n,t+1}-T_{n,t}})=\\(P_{n,t}^{\text{Heat}}-P_{n,t}^{\text{Diss}}-P_{n,t}^{\text{Cool}})\Delta t \label{eq:thermal}
		\end{split}\end{equation}\end{minipage}\\
		& \begin{minipage}{0.4\textwidth}\begin{equation}\text{Electrolytic heat: }P_{n,t}^{\text{Heat}}=P_{n,t}^{\text{Stack}}-N^{\text{Cell}}U^{\text{tn}}I_{n,t}  \label{eq:Pgen}
		\end{equation}\end{minipage}\\
		& \begin{minipage}{0.4\textwidth}\begin{equation}\text{Heat dissipation: }P_{n,t}^{\text{Diss}}=({T_{n,t}-T^{\text{Am}}})/{R_n^{\text{Diss}}}
		\end{equation}\end{minipage}\\
		& \begin{minipage}{0.4\textwidth}\begin{equation}\text{Active cooling: }P_{n,t}^{\text{Cool}} \leq (1-b^{\text{Idle}}) c^{\text{Cool}}(T_{n,t}-T^{\text{Cool}}) \label{eq:Pcool}
		\end{equation}\end{minipage}\\
		& \begin{minipage}{0.4\textwidth}\begin{equation}\text{Temperature limits: }\underline T\leq T_{n,t} \leq \overline T  \label{eq:T}
		\end{equation}\end{minipage}\\
		\vspace{-9pt}
		\\
		\hline
		\multirow{4}{*}{\tabincell{c}{State\\switching}}
		& \begin{minipage}{0.4\textwidth}\begin{equation}\text{Logical constraint: }b_{n,t}^{\text{On}}+b_{n,t}^{\text{By}}+b_{n,t}^{\text{Idle}}=1 \label{eq:logic}
		\end{equation}\end{minipage}\\
		& \begin{minipage}{0.4\textwidth}\begin{equation}\text{Startup action: }b_{n,t}^{\text{On}}+b_{n,t}^{\text{By}}+b_{n,t-1}^{\text{Idle}}-1\leq b_{n,t}^{\text{SU}} \label{eq:SU}
		\end{equation}\end{minipage}\\
		& \begin{minipage}{0.4\textwidth}\begin{equation}\text{Shutdown action: }b_{n,t-1}^{\text{On}}+b_{n,t-1}^{\text{By}}+b_{n,t}^{\text{Idle}}-1\le b_{n,t}^{\text{SD}} \label{eq:SD}
		\end{equation}\end{minipage}\\
		& \begin{minipage}{0.4\textwidth}\begin{equation}\text{Startup delay: }-b_{n,t-2}^{\text{Idle}}-b_{n,t}^{\text{Idle}}+b_{n,t-1}^{\text{Idle}} \leq 0 \label{eq:Idle}
		\end{equation}\end{minipage}\\
		\vspace{-9pt}
		\\
		\hline\hline
	\end{tabular}\vspace{2pt}
	\raggedright \footnotesize Note:
	$U^{\text{rev}}$ and $U^{\text{tn}}$ are the reversible voltage and thermal neutral voltage; $N^{\text{cell}}$ and $A$ are the cell number and electrode area of the stack; $r_1,r_2,s_1,$ $t_1,t_2$ and $t_3 $ are coefficients of the U-I curve; $F$ is the Faraday constant; $f_1$ and $f_2$ are the coefficients of Faradaic efficiency.
	\vspace{-6pt}
\end{table}

\subsubsection{State switching submodel}

ELZs switch among \emph{production}, \emph{standby}, and \emph{idle} states to achieve a wide and flexible load range to accommodate the volatile renewable power, similar to power system unit commitment (UC).
Their transitions are governed by (\ref{eq:logic})--(\ref{eq:Idle}) \cite{varela2021modeling}.  For detailed ELZ models, interested readers can refer to our prior work \cite{qiu2023extend,zeng2024scheduling}.

\subsection{Reactive Power Characteristics}
\label{sec:reactpower}

\begin{figure}[t]
	\centering
	\includegraphics[width=3.45in]{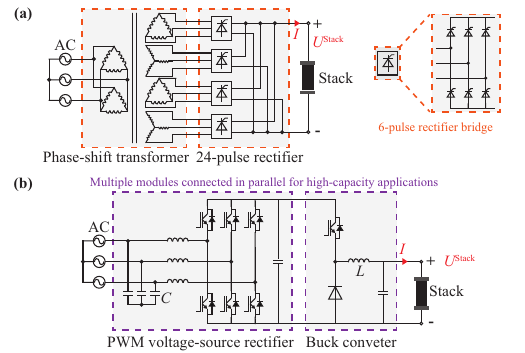}\vspace{-0pt}
	\caption{The detailed topologies of (a) 24-pulse TR, and (b) IGBT-R.}
	\label{fig:TRtopology}\vspace{-0pt}
\end{figure}

We now compare ELZs powered by TRs and IGBT-Rs in terms of reactive power and efficiency. Fig. \ref{fig:TRtopology} shows the topologies of a 24-pulse TR and an IGBT-R.

\subsubsection{TR}

The TR consumes a significant amount of reactive power, which is coupled with the active power and exhibits a nonlinear relationship. The relation between the reactive power $Q^\text{TR}$, stack voltage $U^{\text{Stack}}$ and AC bus voltage $U^{\text{AC}}$ follows \cite{zeng2024scheduling}:
\begin{align}
	Q^\text{TR} &= 2.44{P^\text{Stack}U^\text{AC}}/(\eta^{\text{Rec}}{K^{\text{Rec}}U^\text{Stack}}) \nonumber\\
	&\times\sqrt{\sin^2\left[\arccos(\frac{K^{\text{Rec}}U^\text{Stack}}{2.44U^\text{AC}})\right]+\frac{1-\nu^2}{\nu^2}},  \label{eq:Qscr}
\end{align}
where $K^{\text{Rec}}$ is the turn ratio of the rectifier transformer and where $\nu$ is the harmonic factor of the 24-pulse TR.
{\color{black}The detailed derivation can be found in \cite{zeng2024scheduling}.}

{\color{black}
	Since $P^\text{Stack}$ and $U^\text{Stack}$ are functions of $I$ and $T$, substituting (\ref{eq:Pstack}) and (\ref{eq:Ustack}) into (\ref{eq:Qscr}) yields the complete expression for $Q^{\text{TR}}(U^{\text{AC}}, I, T)$, which is employed for MR configuration planning in Section \ref{sec:SOD}.
	As the OLTC can regulate the AC voltage $U^\text{AC}$ at the hydrogen plant, thereby affecting the TR's firing angle and its reactive power demand, coordinated control of the OLTC and TR can be used to optimize system-level operation. The detailed formulation is presented in Table \ref{tab:scheduling} and discussed in Section \ref{sec:SOD}.
	
}

\subsubsection{IGBT-R}
\label{sec:IGBT-R}

Compared with TRs, IGBT-Rs offer more flexible reactive power characteristics. Two operational modes exist 
as listed below.
\begin{itemize}
	\item \emph{Fixed Power Factor (PF = 1)} (Industrial practices \cite{koponen2021comparison}): Many existing distributed control systems (DCSs) of hydrogen plants lack a built-in reactive power control interface. Thus, IGBT-Rs are often set to a fixed unit PF: 
	\begin{equation}
		Q^{\text{IGBT}}=0. \label{eq:Qigbt1}
	\end{equation}
	\item \emph{Adjustable Reactive Power} (Academia \cite{de2023hydrogen,zhang2022research} and industry standards \cite{biaozhun}): IGBT-Rs can supply or absorb reactive power within capacity limits \cite{de2023hydrogen} to act as a compensator: 
	\begin{equation}
		(P^{\text{Rec}})^2+(Q^{\text{IGBT}})^2\leq(S^\text{IGBT})^2. \label{eq:Qigbt2}
	\end{equation}
\end{itemize}

\subsection{{\color{black}Performance Comparison of Different MR Configurations}}
\label{sec:comparision}

{\color{black}This section compares different examples to show that intuitive rectifier configurations may fall short, while an MR setup offers better performance.}

\begin{figure}[t]
	\includegraphics[width=3.35in]{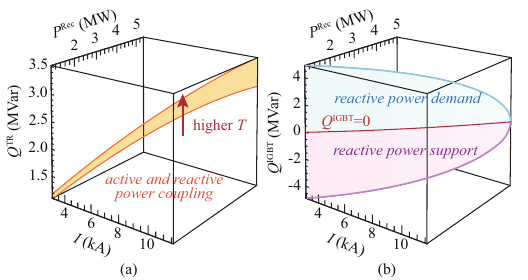}\vspace{-0pt}
	\caption{The active-reactive power characteristics of ELZs powered by different rectifiers under varying electrolytic currents. (a) TR. (b) IGBT-R}
	\label{fig:region}
	\vspace{-0pt}
\end{figure}

\begin{figure}[t]   \centering
	\includegraphics[width=3.3in]{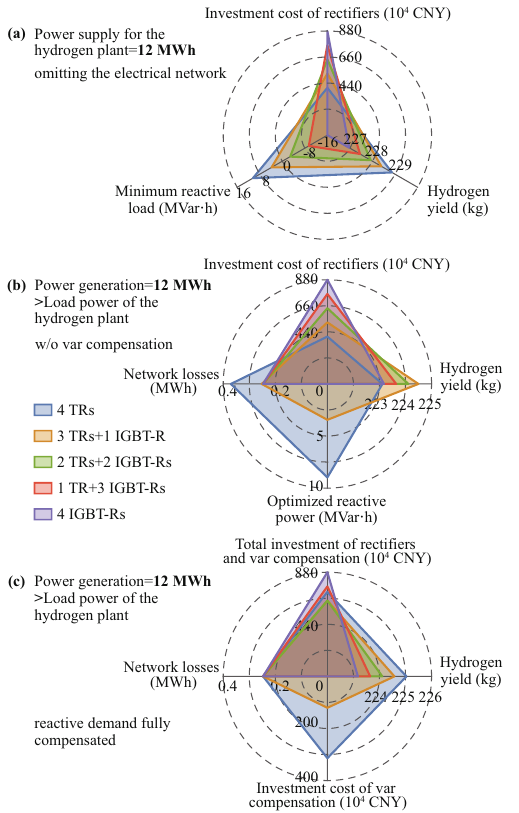}\vspace{-0pt}
	\caption{{\color{black} Performance comparison of different rectifier configurations for four ELZs. (a) A hydrogen plant without considering the network constraints. (b) A simple illustrative system without var compensation. (c) A simple illustrative system with the reactive demand fully compensated.}}
	\label{fig:comparision}
	\vspace{-0pt}
\end{figure}

\subsubsection{{\color{black}A single rectifier}}

\textcircled{1} As shown in Fig. \ref{fig:TRtopology}, the PWM and two-stage topology of the IGBT-R result in higher switching and conduction losses than those of the TR, leading to greater losses $P^\text{Loss}$ and lower efficiency $\eta^{\text{Rec}}$ than those of the TR. \textcircled{2} The IGBT-R also requires a larger plant area, increasing investment costs. \textcircled{3} Figs. \ref{fig:region}(a) and \ref{fig:region}(b) compare the active and reactive power characteristics of the TR and IGBT-R for a 5 MW ELZ (rated flow of 1,000 Nm$^3$/h). As shown in Fig. \ref{fig:region}(a), the TR's active and reactive power are coupled, meaning that for a fixed electrolytic current, the reactive power remains uncontrollable, varying only within the orange region as the stack temperature changes. In contrast, Fig. \ref{fig:region}(b) shows that the IGBT-R provides greater flexibility. For a fixed active power or current, the IGBT-R can supply reactive power (purple region), absorb reactive power (blue region), or maintain a PF of 1 (red curve).

\begin{figure}[t]
	\centering
	\includegraphics[width=3.35in]{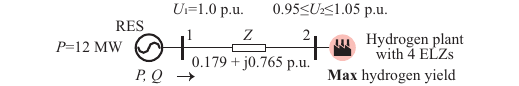}\vspace{-0pt}
	\caption{{\color{black}A simple illustrative system for comparing the performance of different rectifier configurations.}}
	\label{fig:example}\vspace{-0pt}
\end{figure}

{\color{black}
	\subsubsection{MR configuration}

	To illustrate the impact of different rectifier configurations, three examples are given.
	In the first example, a hydrogen plant with four 5 MW ELZs supplied with 12 MW of power for one hour. Using hydrogen production maximization as the objective, we calculate the hydrogen yield and reactive demand for each configuration. Fig. \ref{fig:comparision}(a) presents the results along with the investment costs. TRs offer lower costs and slightly higher hydrogen yield. If the interaction between hydrogen production and the electrical network is omitted, the all-TR setup seems the optimal choice. However, as discussed in Section \ref{sec:ReP2H}, both voltage security and network losses must be carefully addressed. Thus, further analysis is carried out based on a simple illustrative system, as shown in Fig. \ref{fig:example}.
	
	In the second example, the power source remains at 12 MW. However, due to network losses, the actual power delivered to the hydrogen plant is reduced. The optimized performance under different configurations is presented in Fig. \ref{fig:comparision}(b). In contrast to the first case, the all-TR setup results in the lowest hydrogen production, primarily due to significant network losses induced by high reactive power demand. Since IGBT-Rs are capable of providing reactive power support, incorporating them enhances overall system performance. When considering both investment cost and hydrogen yield, the 4-TR setup is no longer optimal. An MR configuration improves hydrogen yield with a modest cost increase, and the 3 TR + 1 IGBT-R configuration achieves the best performance.
	
	In the third example, voltage security and network losses are explicitly considered, making var compensation essential. To isolate its impact, we assume full compensation for the hydrogen plant's reactive load. The corresponding results are shown in Fig. \ref{fig:comparision}(c). Although the 4-TR configuration delivers the highest productivity, the additional cost of the var compensator increases the total investment. In comparison, using one IGBT-R achieves similarly good productivity with lower total cost. Together with the second example, this confirms the necessity of coordinating hydrogen production with the electrical network, offering insights that challenge intuitive preferences. 
	
	It should be noted, however, that these examples represent only a single-time snapshot. Whether an MR configuration can optimally balance these factors over a longer continuous operation period requires further investigation.

}


\section{Optimal Investment Portfolio Planning of Mixed Rectifier Configurations}
\label{sec:SOD}


\subsection{Planning {\color{black}and Operational} Framework}

{\color{black}The optimal investment portfolio of MRs is a mixed planning and operation problem, consisting of two stages and their interaction}:
\begin{itemize}
  \item \emph{Investment stage}: Determines the number of ELZs, type and number of rectifiers, and capacity of the static var generator (SVG) for compensation.
  \item \emph{Operational stage}: Employs the coordinated active-reactive power scheduling method proposed in our prior work \cite{zeng2024scheduling} to manage renewables, multiple ELZs, BES, and var resources.
\end{itemize}

\begin{figure}[t]
  \centering
  \includegraphics[width=3.45in]{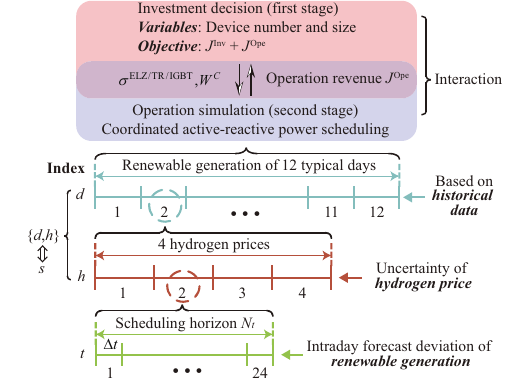}\vspace{-0pt}
  \caption{{\color{black}Planning and operational framework for mixed rectifiers, ELZs, and var resources in the ReP2H system.}}   \label{fig:LST}\vspace{-0pt}
\end{figure}

{\color{black}
In the framework summarized in Fig. \ref{fig:LST}, the investment decisions are passed to the operation stage once determined. The operation stage then performs simulations with the objective of maximizing revenue and feeds the result back to the investment stage, similar to the procedure adopted in \cite{jiang2025two}.}

The ReP2H system faces uncertainties \cite{li2022coordinated}, including the hydrogen price and renewable power over three different timescales. We address these issues using the framework in Fig. \ref{fig:LST}. First, we use 12 typical days across four seasons to represent baseline patterns of renewable power. Second, for each typical day, we use four hydrogen price scenarios \cite{yang2024stochastic} to describe its uncertainty. 
Finally, intraday wind and solar power uncertainties are handled via IGDT \cite{yu2024optimal}. The intraday uncertainty is modeled as a set $\mathcal{B}$, which describes the uncertainty horizon $\alpha $ on the basis of the baseline $\tilde{P}^{{\text{WT/PV}}}_{s,t}$, formulated as: \begin{align}
    \mathcal{B}(\alpha,\tilde{P}_{s,t}^{{\text{WT/PV}}})=\Big\{P_{s,t}^{{\text{WT/PV}}}:\Big|\frac{P_{s,t}^{{\text{WT/PV}}}-\tilde{P}_{s,t}^{{\text{WT/PV}}}}{\tilde{P}_{s,t}^{{\text{WT/PV}}}}\Big|\leq \alpha \Big\}.
\end{align} where $s$ represents the day price index $\{d,h\}$ for simplicity.

\subsection{Objective Function}

The objective of the MR investment portfolio model is to maximize the total revenue $J^{\text{Tot}}$:
\begin{align}
&\max~J^{\text{Tot}}= J^{\text{Inv}}+J^{\text{Ope}}, \label{eq:Ctot}
\end{align}
where the investment cost $J^{\text{Inv}}$ covers ELZs, rectifiers, and var resources, as
\begin{align}
\nonumber J^{\text{Inv}}&=-{r(1+r)^y}/[{(1+r)^y-1}]  \label{eq:Cinv} \\
  &\times (c^{\text{C}}W^{\text{C}}+c^{\text{TR}}\sigma^{\text{TR}}+c^{\text{IGBT}}\sigma^{\text{IGBT}}+c^{\text{ELZ}}\sigma^{\text{ELZ}}),
\end{align}
and the operational revenue $J^{\text{Ope}}$ considers hydrogen production and sale revenues, startup/shutdown costs of ELZs, and operation and maintenance (O\&M) costs, i.e.,
\begin{align}
    J^{\text{Ope}}=&365\sum_{s=1}^{N_{s}}\pi_{s} J_{s}^{\text{Ope}}, \label{eq:Cope}\\
   J_{s}^{\text{Ope}}=&\sum_{t=1}^{N_t}\sum_{n=1}^{\sigma^{\text{ELZ}}}(c^{\text{H}_2}Y_{n,t}^{\text{H}_2}-c^{\text{SU}}b_{n,t}^{\text{SU}}-c^{\text{SD}}b_{n,t}^{\text{SD}})
+\frac{\eta_\text{O\&M}J^{\text{Inv}}}{365}.    \label{eq:Csope}
\end{align}
Here, $r$ and $y$ are the discount rate and equipment lifespan, respectively; $\pi_{s}$ represents the probability of scenario $s$; and $\eta_\text{O\&M}=2\%$ is the annual O\&M cost factor. {\color{black}In (\ref{eq:Csope}), maximizing hydrogen yield simultaneously achieves maximum RE utilization and reactive power optimization, while also accounting for SVG planning and network losses.}

\subsection{Investment and Operational Constraints}

\subsubsection{Investment constraints}

Constrained by available land, the number of ELZs satisfies (\ref{eq:Nelz}), and the number of rectifiers must match the ELZs, following (\ref{eq:Nrec}). To ensure the electrical strength of the AC system, the short-circuit ratio $R^{\text{sc}}$ at the hydrogen plant must exceed a threshold \cite{liu2024system}. TRs rely on grid support for reliable commutation, requiring high $R^{\text{sc}}$, whereas IGBT-Rs can provide support to the AC grid \cite{zhang2022research,tavakoli2023gridforming}, reducing the requirement for $R^{\text{sc}}$. Thus, $R^{\text{sc}}$ is estimated as (\ref{eq:scr}) for a fast assessment of electrical strength. 
\begin{align}
  &0\leq \sigma^{\text{ELZ}}\leq \overline{\sigma}^{\text{ELZ}}, \label{eq:Nelz}\\
  &\sigma^\text{TR}+\sigma^\text{IGBT}=\sigma^{\text{ELZ}}, \label{eq:Nrec}\\
  &R^{\text{sc}}=\frac {S^{\text{sc}}}{(\sigma^\text{TR}S^\text{TR}+\sigma^\text{IGBT}S^\text{IGBT})}\geq \frac {3\sigma^\text{TR}+2\sigma^\text{IGBT}}{\sigma^\text{TR}+\sigma^\text{IGBT}}. \label{eq:scr}
\end{align}

Note that the electrical strength is roughly evaluated here for investment planning. An accurate assessment of the strength and threshold of the short-circuit ratio in ReP2H systems requires electromagnetic analysis and should be explored in future work.

\begin{table*}[tb]\scriptsize
  \renewcommand{\arraystretch}{1.58}
  \caption{Model and Operational Constraints of the AC Network in the ReP2H System}\vspace{-0pt}
  \label{tab:scheduling}
  \centering
  \begin{tabular}{c@{\hspace{3pt}}c@{\hspace{9pt}}c}
  \hline\hline
            & Operational constraints                  \\
  \hline
  \multirow{8}{*}{Power flow \cite{farivar2013branch}}
  &\begin{minipage}{0.84\textwidth}  \begin{equation}\text{Branch active power injections: }\textstyle\sum_{j^{\prime}\in\sigma(j)}P_{jj^{\prime},t}=p_{j,t}+\sum_{i\in\pi(j)}P_{ij,t}-r_{ij}\ell_{ij,t} \label{eq:branchflow}
   \end{equation}\end{minipage}\\
    & \begin{minipage}{0.84\textwidth}\begin{equation}\text{Branch reactive power injections: } \textstyle\sum_{j^{\prime}\in\sigma(j)}Q_{jj^{\prime},t}=q_{j,t}+\sum_{i\in\pi(j)}Q_{ij,t}-x_{ij}\ell_{ij,t}
   \end{equation}\end{minipage}\\
   &\begin{minipage}{0.84\textwidth}\begin{equation}\text{Ohm's law of transmission line: }\upsilon_{j,t}=\upsilon_{i,t}-2(r_{ij}P_{ij,t}+x_{ij}Q_{ij,t})+(r_{ij}^2+x_{ij}^2)\ell_{ij,t}
   \end{equation}\end{minipage}\\
   & \begin{minipage}{0.84\textwidth}\begin{equation}\text{Ohm's law of transformer branch: } \textstyle\sum_{k=0}^{K_{ij}}(w_{ij,k})^2\delta_{ij,k,t}\upsilon_{j,t}=\upsilon_{i,t} -2(r_{ij}P_{ij,t}+x_{ij}Q_{ij,t})+(r_{ij}^2+x_{ij}^2)\ell_{ij,t} \end{equation}\end{minipage}\\
   &\begin{minipage}{0.84\textwidth}\begin{equation}\text{Second-order cone relaxation: }\left\|(2P_{ij,t},2Q_{ij,t},\ell_{ij,t}-\upsilon_{i,t})\right\|_2\leq\ell_{ij,t}+\upsilon_{i,t} \end{equation}\end{minipage}\\
  &\begin{minipage}{0.84\textwidth}\begin{equation}\text{Auxiliary variables for squared voltage and current and security constraints: }\underline{U}_j^2\leq\upsilon_{j,t}=\left|U_{j,t}\right|^2\leq\overline{U}_j^2,\ell_{ij,t}=\left|I_{ij,t}\right|^2\leq\overline{I}_{ij}^2
   \end{equation}\end{minipage}\\
   &\begin{minipage}{0.84\textwidth}\begin{equation}\text{Bus injected power: }p_{j,t}=P_{j,t}^{\text{WT}}+P_{j,t}^{\text{PV}}-P_{j,t}^{\text{BES,C}}+P_{j,t}^{\text{BES,D}}-p_{j,t}^{\text{L}},~
  q_{j,t}=Q_{j,t}^\text{WT}+Q_{j,t}^\text{PV}+Q_{j,t}^\text{BES}+Q_{j,t}^\text{C}-q_{j,t}^\text{L}  \label{eq:powerinjection}
   \end{equation}\end{minipage}\\
   &\begin{minipage}{0.84\textwidth}\begin{equation}\text{Power of hydrogen plant: }p_{j,t}^{\text{L}}=\textstyle\sum_{n=1}^{{\sigma^{\text{ELZ}}}} P_{n,t}^{\text{ELZ}}, ~q_{j,t}^\text{L}=\textstyle\sum_{n=1}^{{\sigma^{\text{TR}}}}Q_{n,t}^\text{TR}+\sum_{n=1}^{{\sigma^{\text{IGBT}}}}Q_{n,t}^\text{IGBT} \label{eq:H2plant}
   \end{equation}\end{minipage} \\
   \vspace{-10.5pt}
   \\
   \hline
   \multirow{2}{*}{\tabincell{c}{WTs and\\PV plants \cite{zeng2024scheduling}}}
   & \begin{minipage}{0.84\textwidth}\begin{equation}\text{Reactive power of WTs: }1.24P_{j,t}^\text{WT} - 0.91S_{j}^\text{WT}\le Q_{j,t}^\text{WT} \le0.91S_{j}^\text{WT}-0.58P_{j,t}^\text{WT} \label{eq:qwt}
   \end{equation}\end{minipage}\\
   & \begin{minipage}{0.84\textwidth}\begin{equation}\text{Reactive power of PV plants: }
    (P_{j,t}^{\text{PV}})^2 + (Q_{j,t}^{\text{PV}})^2\leq (S_j^{\text{PV}})^2,~
    -P_{j,t}^{\text{PV}}\tan\theta\le Q_{j,t}^{\text{PV}}\le P_{j,t}^{\text{PV}}\tan\theta  \label{eq:qpv}
   \end{equation}\end{minipage}  \\
   \vspace{-10.5pt}
   \\
   \hline
   \multirow{5}{*}{BES \cite{ghazavi2021simultaneous}}
   & \begin{minipage}{0.84\textwidth}\begin{equation}\text{Active and reactive power of BES: }(P_{j,t}^{\text{BES,C/D}})^2+(Q_{j,t}^{\text{BES}})^2\le(S_j^{\text{BES}})^2 \label{eq:qbes}
   \end{equation}\end{minipage}\\
   & \begin{minipage}{0.84\textwidth}\begin{equation}\text{Charge and discharge: }0\le P_{j,t}^{\text{BES,C/D}}\le b_{j,t}^{\text{BES,C/D}}\overline{P}^{\text{BES,C/D}}, ~
  b_{j,t}^{\text{BES,C}}+b_{j,t}^{\text{BES,D}}\le 1  \label{discharge}
   \end{equation}\end{minipage}\\
   & \begin{minipage}{0.84\textwidth}\begin{equation}\text{State of charge of BES: }E_{j,t+1}=(1-\zeta^\text{BES})E_{j,t}+(\eta^{\text{BES,C}}P_{j,t}^{\text{BES,C}}- {P_{j,t}^{\text{BES,D}}}/{\eta^{\text{BES,D}}}){\Delta t} \label{soc1}
   \end{equation}\end{minipage}\\
   & \begin{minipage}{0.84\textwidth}\begin{equation}\text{State of charge constraints: }E_{j, t=0}=E_{j,{t=N_t}},~
  \underline{E}\leq E_{j,t}\leq\overline{E} \label{soc2} \vspace{0.5pt}
   \end{equation}\end{minipage}\\
   \vspace{-10.5pt}
   \\
   \hline
   Var compensation
   & \begin{minipage}{0.84\textwidth}\begin{equation}\text{Reactive power of var compensation: }-W^{\text{C}}\le Q_{j,t}^\text{C}\le W^{\text{C}} \label{Qc} \end{equation}\end{minipage}\\
    \vspace{-10.5pt}
   \\
   \hline
   \multirow{2}{*}{OLTC \cite{zeng2024scheduling}}
   & \begin{minipage}{0.84\textwidth}\begin{equation}\text{OLTC tap position modeling: }k_{ij,t}=\textstyle\sum_{k=0}^{K_{ij}}\delta_{ij,k,t}w_{ij,k},~
   \sum_{k=0}^{K_{ij}}\delta_{ij,k,t}=1,~\delta_{ij,k,t}\in\{0,1\}  \label{eq:oltc1}
   \end{equation}\end{minipage}\\
   & \begin{minipage}{0.84\textwidth}\begin{equation}\text{Action and operational range constraints: }\textstyle\sum_{t=2}^{{N_t}}\left|k_{ij,t}-k_{ij,t-1}\right|\le\overline{k}_{ij}^{\text{All}},~
\underline{k}_{ij}\le k_{ij,t}\le\overline{k}_{ij} \label{eq:oltc2}
   \end{equation}\end{minipage} \\
   \vspace{-10.5pt}
   \\
  \hline\hline
  \end{tabular}\vspace{-0pt}
\end{table*}
\subsubsection{Operational constraints}

We consider the coordination among active and reactive resources such as the electrical network, ELZs, WTs, PV plants, BES, var compensation, and OLTC. Their models are presented in Table \ref{tab:scheduling}, where we use the DistFlow model \cite{farivar2013branch} to describe the network power flow, which is radial in common ReP2H systems.

In operation, active power balance is achieved by on-off switching and load allocation among multiple ELZs in the hydrogen plant, assisted by the BES. The reactive power flow is optimized in coordination with the active power to support the voltage and reduce network losses. 
For details of the operation model, readers are referred to \cite{zeng2024scheduling}.

\subsection{The Overall Optimal Portfolio Model}
\label{sec:igdt}

\subsubsection{Planning with $\alpha =0$}

As per the IGDT procedure, first, we assume zero intraday uncertainty for the wind and solar power in typical-day scenarios. The optimal portfolio model for the MR configuration is given as follows:
\begin{subequations}\label{eq:plan1}
	\begin{align}
		& \tilde{J}^{\text{Tot}} = \max~\text{(\ref{eq:Ctot})}, \\
		& \text{s.t.} ~\text{(\ref{eq:Pelz})-(\ref{eq:Idle}), (\ref{eq:Qscr}), (\ref{eq:Qigbt1}) or (\ref{eq:Qigbt2}),
			(\ref{eq:Nelz})-(\ref{eq:scr}),
			(\ref{eq:branchflow})-(\ref{eq:oltc2}),} \label{eq:conserplancons}
	\end{align}
\end{subequations}
where $\tilde{J}^{\text{Tot}}$ is the total revenue. We denote the investment cost as $\tilde{J}^{\text{Inv}}$ and the operational revenue under scenario $s$ as $\tilde{J}_s^{\text{Ope}}$.

\subsubsection{W-IGDT-based portfolio model}

{\color{black}
	Although the investment decisions deal with a long-term time frame considering the lifespan of equipment, the optimization problem incorporates a short-term scheduling program, in which the operator needs to handle the intraday uncertainty of renewable power. To address the short-term uncertainty,} a risk-averse model based on robust IGDT \cite{yu2024optimal} is developed.
Based on the revenue of (\ref{eq:plan1}), with a preset deviation factor $\beta>0$, the robust IGDT model maximizes the uncertainty boundary $\alpha$ under the worst-case renewable generation, as follows:
\begin{subequations}\label{eq:IGDT}
	\begin{align}
		\max~&\alpha \\
		\text{s.t.}~~&J^{\text{Tot}}\geq(1-\beta)\tilde{J}^{\text{Tot}},\\
		&P_{s,t}^{{\text{WT/PV}}}=(1-\alpha)\tilde{P}_{s,t}^{{\text{WT/PV}}},~\forall s\\
		&\text{(\ref{eq:conserplancons})}.
	\end{align}
\end{subequations}

Since higher renewable penetration corresponds to higher risk, scenarios with higher renewable power output require larger uncertainty boundaries to mitigate the risk \cite{nasr2020assessing}. Therefore, we assign weights to the boundaries $\alpha_s$ on the basis of power output levels for different scenarios, improving the conventional IGDT model to the W-IGDT \cite{nasr2020assessing} model to better manage the risk. The W-IGDT version of (\ref{eq:IGDT}) is as follows:
\begin{subequations}\label{eq:W-IGDT}
	\begin{align}
		&\max~ \sum_{s}\pi_s\frac{N_{s}\sum_t( \tilde{P}_{s,t}^{{\text{WT}}}+\tilde{P}_{s,t}^{{\text{PV}}})}{\sum_{s}\sum_t(\tilde{P}_{s,t}^{{\text{WT}}}+\tilde{P}_{s,t}^{{\text{PV}}})}\alpha_s, \\
		&~\text{s.t.}~~J^{\text{Inv}}+365J_s^{\text{Ope}}\geq(1-\beta)(\tilde{J}^{\text{Inv}}+365\tilde{J}_s^{\text{Ope}}),~\forall s,\\
		&~~~~~~P_{s,t}^{{\text{WT/PV}}}=(1-\alpha_s)\tilde{P}_{s,t}^{{\text{WT/PV}}},~\forall s, \\
		&~~~~~~\text{(\ref{eq:conserplancons}).}
	\end{align}
\end{subequations}


\subsection{Two-stage Decomposition-Based Solution Approach---PH}
\label{sec:PH}

The established optimal portfolio model (\ref{eq:plan1}) uses mixed-integer nonlinear programming (MINLP), which is difficult to solve. Therefore, we use polynomial approximation, piecewise linearization, and the Big-M method to transform it into mixed-integer second-order cone programming (MISOCP). Nevertheless, the reformulated model is large in scale, involves various active and reactive power resources, and requires multiple scenarios to address uncertainty. Specifically, for an 8-node small-scale system with 4 ELZs (in Section \ref{sec:Systems}-1), the model has 358,512 constraints, 167,328 continuous variables, and 42,624 binary variables. This makes direct solving infeasible, akin to large-scale systems.

Fortunately, the proposed model is a standard two-stage SP, with each stage presenting investment or operation. Therefore, we can adopt the PH algorithm \cite{rockafellar1991scenarios}, which is effective for accelerating the solution of the SP. We decompose the SP into $N_s$ subproblems, compactly expressed as:
\begin{align}
	\min_{\{\bm{x},\bm{y}_s,\forall s\}}&\bm{a}^\text{T}\bm{x}+\sum_{s=1}^{N_s}\pi_s \bm{b}_s^\text{T}\bm{y}_s  \label{eq:PHobj}\\
	\text{s.t.}~~~&\bm{x}\in X, \bm{y}_s\in Y_s, \forall s  \label{eq:PHst}
\end{align}
where $\bm{a}$ and $\bm{b}_s$ are coefficient matrices; $\bm{x}$ and $\bm{y}_s$ are the vectors of investment and operational variables, respectively, with $X$ and $Y_s$ as their feasible sets. Then, we can use Algorithm \ref{alg:PH} to solve the portfolio model effectively by solving the investment and operational variables iteratively.

\begin{figure}[t]
  \centering
  \includegraphics[width=3.45in]{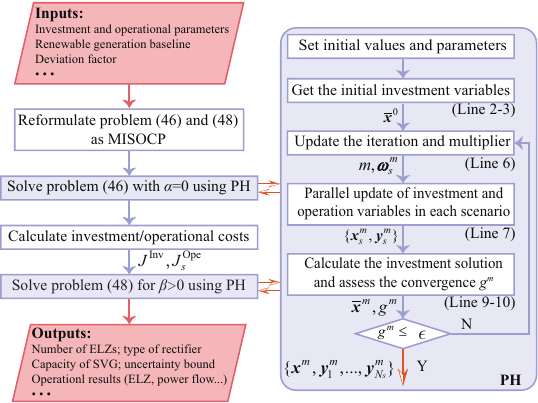}\vspace{-0pt}
\caption{\color{black}Solution procedure of the proposed method for MR configuration and var resource planning.}
  \label{fig:procedure}\vspace{-0pt}
\end{figure}

{\color{black}
	The overall solution procedure of the proposed method and its application to MR configuration and var resource planning in the ReP2H system are summarized in Fig. \ref{fig:procedure}.
}





\begin{algorithm}[t]
  \caption{PH Algorithm for the Portfolio Model}
  \label{alg:PH}
  \begin{algorithmic}[1]
    \STATE \textbf{precondition:} Initialize the iteration index $m\gets0$, Lagrange multiplier $\bm{\omega}_s^0\gets0$ and convergence index $g^0$, and select the convergence tolerance $\epsilon \gets10^{-3}$ \label{line:1}
    \STATE  Get the investment and operational variables in parallel by   $\{\bm{x}_s^0,\bm{y}_s^0\}\leftarrow\arg\min_{\{\bm{x},\bm{y}_s\}}\{\bm{a}^\text{T}\bm{x}+\bm{b}_s^\text{T}\bm{y}_s:
  \text{(\ref{eq:PHst})}\}$ for $\forall s$ \label{line:3}
    \STATE Calculate the initial value for the investment variable $\bm{\overline{x}}^0\leftarrow\textstyle\sum_s\pi_s\bm{x}_s^0$ \label{line:5}
    \WHILE{$g^m\leq\epsilon$}{
    \FOR {$s=1$ \TO $N_s$}
    \STATE Update the index $m\gets m+1$, and multiplier $\bm{\omega}_s^m\leftarrow\bm{\omega}_s^{m-1}+\rho(\bm{x}_s^{m-1}-\bm{x}^{m-1})$ \label{line:8}
    \STATE Update the investment and operational variables $\{\bm{x}_s^m,\bm{y}_s^m\}\gets\arg\min_{\{\bm{x},\bm{y}_s\}}\{\bm{a}^\text{T}\bm{x}+\bm{b}_s^\text{T}\bm{y}_s+(\bm{\omega}_s^m)^\text{T}\bm{x}
+\frac{\rho}{2}\|\bm{x}-\overline{\bm{x}}^{m-1}\|_2^2:\text{(\ref{eq:PHst})}\}$ \label{line:9}
    \ENDFOR
    \STATE Calculate the investment solution $\bm{\overline{x}}^m\leftarrow\textstyle\sum_s\pi_s\bm{x}_s^m$  \label{line:11}
    \STATE Assess the convergence $g^m\leftarrow\textstyle \sum_s\pi_s\parallel \bm{x}_s^m-\bm{\overline{x}}^m \parallel_2$ \label{line:12}
    }\ENDWHILE
    \RETURN $\{\bm{\overline{x}}^m,\bm{y}_1^m,\ldots,\bm{y}_{N_s}^m\}$
  \end{algorithmic}
\end{algorithm}


\section{Case Studies}
\label{sec:cases}

\begin{figure}[t]
  \centering
  \includegraphics[width=3.4in]{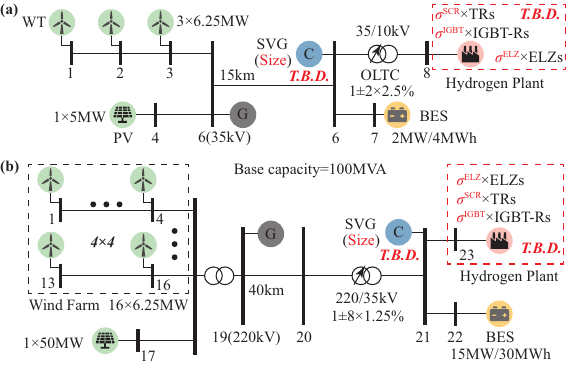}\vspace{-0pt}
  \caption{Topologies of the ReP2H systems used for the case study. (a) Small-scale system with 23.75 MW of RE. (b) Large-scale system with 150 MW of RE.}   \label{fig:case}\vspace{-0pt}
\end{figure}

Case studies are performed based on two real-world systems in Inner Mongolia, China, as shown in Fig. \ref{fig:case}. The test systems are adapted from our prior work \cite{zeng2024scheduling}. The small-scale system validates the proposed portfolio model and explores the potential of MR configurations, whereas the large-scale system verifies the scalability of the model and provides insights for practice. Simulations are performed via \textit{Wolfram Mathematica 14.0}, with optimizations solved via \textit{Gurobi 11.0.0}.

\subsection{Test System Settings}
\label{sec:Systems}

\subsubsection{Small-scale system}\label{sec:smallcase} The system consists of 3$\times$6.25 MW WTs at buses 1-3, a 5 MW PV plant at bus 4, and a 2 MW/4 MWh BES at bus 7. The investment decisions involve the number of ELZs and rectifiers at the hydrogen plant at the bus and the SVG capacity at bus 6.



\subsubsection{Large-scale system} \label{sec:largecase}

Includes 16$\times$6.25 MW WTs, 1$\times$50 MW PV plant, 1$\times$15 MW/30 MWh BES. The hydrogen plant and SVG are located at buses 21 and 23, respectively.

\begin{figure}[t]   \centering
  \includegraphics[width=3.4in]{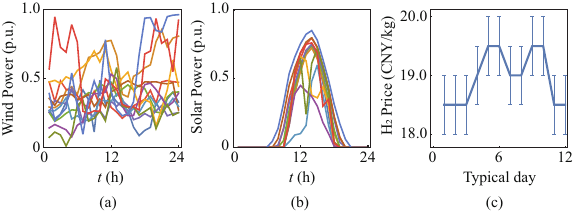}\vspace{-0pt}
  \caption{Scenarios of wind and solar power  and the average hydrogen price of 12 typical days. (a) Wind power. (b) Solar power. (c) Hydrogen price.}
  \label{fig:typicaldays}\vspace{-0pt}
\end{figure}

 \begin{table}[t]\scriptsize
  \renewcommand{\arraystretch}{1.4}
  \caption{The Key Investment Parameters for Case Studies}\vspace{-0pt}
  \label{tab:para}
  \centering
  \begin{tabular}{c@{\hspace{5pt}}c@{\hspace{5pt}}c@{\hspace{5pt}}c@{\hspace{5pt}}c}
  \hline\hline
  Facility                            & Rating       & Unit investment cost         & Discount rate $r$      & Lifetime $y$  \\
  \hline
  ELZ w/o rectifier                 & 5 MW        & 7,000,000 CNY {\color{black}\cite{shuidian2024}}           &  \multirow{4}{*}{8\%}    & \multirow{4}{*}{20 years}     \\
  TR                                & 6 MVA       & 1,000,000 CNY {\color{black}\cite{gao2024advanced}}          &     &                                   \\
  IGBT-R                            & 6 MVA       & 2,200,000 CNY {\color{black}\cite{gao2024advanced}}          &     &                               \\
  SVG                               &             & 350 CNY/kVar  {\color{black}\cite{davoodi2025scalable}}        &     &            \\
  \hline\hline
  \end{tabular}
\end{table}

\begin{table}[t]\scriptsize
  \renewcommand{\arraystretch}{1.4}
  \caption{\color{black}The Operation-related Parameters for Case Studies \cite{zeng2024scheduling}}
  \label{tab:para2}
  \centering \color{black}
  \begin{tabular}{c@{\hspace{7pt}}c@{\hspace{7pt}}c@{\hspace{7pt}}c}
  \hline\hline
  Parameter                             & Value                                 & Parameter                        & Value   \\
  \hline
  $C^{\text{ELZ}}$                      & 7.8$\times10^{7}$J/\textcelsius                  & $R^{\text{Diss}}$                      & 1.06 \textcelsius/kW             \\
  $c^{\text{Cool}}$                     & 17 kW/\textcelsius                               & $T^{\text{Cool}}$ / $\eta^{\text{Cool}}$ & 5 \textcelsius\ /  4.0        \\
  $\overline I$ / $\underline I$        & 12 / 2.4 kA                                      & $\overline T$ / $\underline T$        & 80 / 25 \textcelsius                               \\
   $P^{\text{By}}$                      & 0.5 MW               &$c^{\text{SU}}$ / $c^{\text{SD}}$      & 1000 / 0 CNY \\
  $\overline {E}$ / $\underline {E}$      &0.9 / 0.1    p.u.                          & $E_{j,t=0/{N_t}}$   & 0.5 / 0.5 p.u.                \\
  $\eta^{\text{BES,C/D}}$                & 0.95 / 0.95                              &$\overline{U}_j$ / $\underline{U}_j$   & 1.05 / 0.95  p.u. \\
   $\overline{I}_{ij}$ (35kV)  & 0.37 p.u.                                &$r_{ij}$ / $x_{ij}$ (35kV)                    & 0.0096 / 0.0262 p.u./km     \\
   $\overline{I}_{ij}$ (220kV)  & 1.8 p.u.                                &$r_{ij}$ / $x_{ij}$ (220kV)                  & 0.00025 / 0.000832 p.u./km     \\
  \hline\hline
  \end{tabular}
\end{table}

The key investment parameters {\color{black}and operation-related parameters are listed in Tables \ref{tab:para} and \ref{tab:para2}, respectively. More details on} the topology and operational details of the ReP2H system and ELZs are available in \cite{zeng2024scheduling}. Renewable generation for typical days, derived from on-site historical data via K-means clustering, is shown in Figs. \ref{fig:typicaldays}(a) and \ref{fig:typicaldays}(b).
The hydrogen price scenario is shown in Fig. \ref{fig:typicaldays}(c), with fluctuations modeled based on local agricultural ammonia demand \cite{wu2023multi}.

\subsection{Investment Portfolio Results of the Small-scale System for $\beta=0$ $ (\alpha=0)$ }
\label{sec:planning}

\subsubsection{{\color{black}Investment portfolio results}}

\begin{figure}[t]\centering
  \includegraphics[width=3.35in]{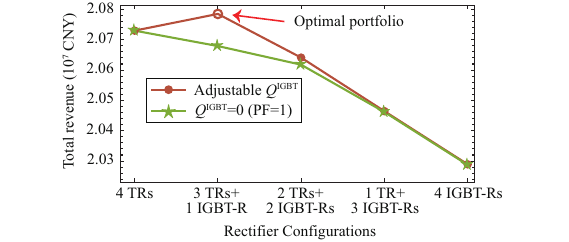}\vspace{-0pt}
  \caption{Total revenue of different rectifier configurations.}   \label{fig:revenue}\vspace{6pt}
  \includegraphics[width=3.35in]{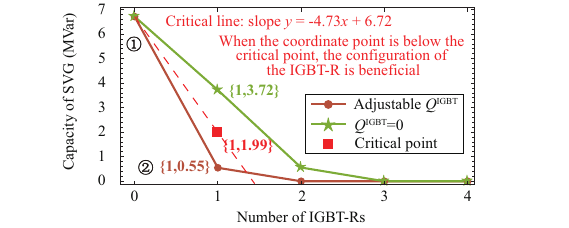}\vspace{-0pt}
  \caption{The optimal SVG capacity under varying numbers of IGBT-Rs.}   \label{fig:SVG}\vspace{-0pt}
\end{figure}

For the small-scale system, with the intraday renewable power uncertainty set to zero, the rectifier configurations of four ELZs are evaluated. Fig. \ref{fig:revenue} presents the total revenues of different portfolios of TR and IGBT-R, while Table \ref{tab:planning} summarizes the planning and operational results. Among the results, cases where IGBT-Rs are allowed to adjust reactive power are denoted \textcircled{2}-\textcircled{5}, and cases where the PF of the IGBT-R is fixed at 1 are denoted \textcircled{6}-\textcircled{9} (see the discussion about reactive power characteristics in Section \ref{sec:reactpower}-2).

As shown in Fig. \ref{fig:revenue} and {\color{black}Table \ref{tab:planning}}, when the IGBT-Rs are allowed to provide reactive power, the optimal configuration is 3 TRs and 1 IGBT-R, {\color{black}which better balances investment and operation}. In contrast, if the IGBT-R is fixed at $\text{PF}=1$, a uniform configuration of 4 TRs is preferable. The reasons are as follows.

Comparing \textcircled{1} and \textcircled{6}, we find that while 1 IGBT-R at $\text{PF}=1$ reduces the need for 4 MVar of SVG, it does not offset its higher cost and lower efficiency, leading to reduced hydrogen yield. Consequently, under configurations \textcircled{6}-\textcircled{9}, IGBT-Rs with $\text{PF}=1$ are less cost-effective than TRs are. By further comparing \textcircled{1} and \textcircled{2}, we find that IGBT-Rs with adjustable reactive power support reduce SVG investment, resulting in lower overall costs for rectifiers and var compensation (\textcircled{1}: 100$\times$4 + 35$\times$6.72 = 635.2$\times$10$^4$ CNY, \textcircled{2}: 100$\times$3 + 220$\times$1 + 35$\times$0.55 = 539.25$\times$10$^4$ CNY), resulting in optimal returns.

In terms of hydrogen production in configurations \textcircled{1} to \textcircled{9}, Table \ref{tab:planning} shows that each additional IGBT-R reduces the hydrogen yield by approximately 2,500 kg/year, resulting in an estimated decrease in revenue of 4.75$\times10^4$ CNY. On the other hand, the annual cost of an IGBT-R is 0.104$\times$(220--100) = 12.48$\times10^4$ CNY greater than that of a TR considering the discount. 
Therefore, an IGBT-R investment is justified only if it reduces SVG capacity by at least (12.48 + 4.75)/(0.104$\times$35) = 4.73 MVar. In particular, Fig. \ref{fig:SVG} shows the SVG capacity with varying numbers of IGBT-Rs, along with the critical point at which configuring more IGBT-Rs is beneficial. It is clear that the only case \textcircled{1}$\rightarrow$\textcircled{2} meets the conditions.

{\color{black}
Moreover, the planning result for the case with 3 ELZs is also discussed, as shown in Table \ref{tab:planning}. It configures 3 TRs to enhance hydrogen production and deploys SVG to meet the reactive power demand. Compared with the optimal 4-ELZ configuration, although the investment is reduced, reduced RE utilization and operational changes lead to lower operational revenue, performing worse than the 4-ELZ setup.
}

\begin{table}[t]\scriptsize
  \renewcommand{\arraystretch}{1.45}
  \caption{Planning and Operational Results under Different Rectifier Configurations in the Small-Scale System}\vspace{-0pt}
  \label{tab:planning}
  \centering
  \begin{threeparttable}
  \begin{tabular}{c@{\hspace{5pt}}c@{\hspace{5pt}}c@{\hspace{5pt}}c@{\hspace{5pt}}c@{\hspace{5pt}}c}
  \hline\hline
          & \tabincell{c}{Configuration\\(TR + IGBT-R)}       & \tabincell{c}{Total Revenue\\{\color{black}$J^{\text{Inv}}+J^{\text{Ope}}$}\\(10$^7$ CNY/yr)}    & \tabincell{c}{SVG\\ Capacity\\(MVar)}    & \tabincell{c}{Hydrogen Yield\\(10$^3$ kg/yr)}   & \tabincell{c}{Network loss\\(MWh/yr)}  \\
  \hline
  \textcircled{1}&4 + 0               & \tabincell{c}{2.0730\\{\color{black}(-0.3569+2.4299)}}       &  6.72         & \textbf{1362.64}  &  1233.1   \\
  \hline
  \textcircled{2}&3 + 1    &  \tabincell{c}{\textbf{2.0785}\\{\color{black}(-0.3469+2.4254)}}         &  0.55         &  1360.27          &  1229.3 \\

  \textcircled{3}&2 + 2    & \tabincell{c}{2.0641\\{\color{black}(-0.3574+2.4215)}}               &  0            &  1358.21          &  1218.0   \\

  \textcircled{4}&1 + 3   &   \tabincell{c}{2.0465\\{\color{black}(-0.3698+2.4163)}}                 &  0            &  1355.46          &  1217.8  \\

  \textcircled{5}&0 + 4    &  \tabincell{c}{2.0289\\{\color{black}(-0.3823+2.4112)}}               &  0            &  1352.77        &  \textbf{1211.8}  \\
  \hline
  \textcircled{6}&3 + 1     & 2.0679                  &  3.72         &  1360.77       &  1238.3 \\

  \textcircled{7}&2 + 2   & 2.0618                  &  0.56         &  1358.09         &  1235.5 \\

  \textcircled{8}&1 + 3    & 2.0462                  &  0            &  1355.32        &  1228.0 \\

  \textcircled{9}&0 + 4    & 2.0288                 &  0          &  1352.68         &   1217.9 \\
  \hline
  &\color{black}3 + 0 (3 ELZs)              & \color{black} \tabincell{c}{2.0578\\(-0.2685+2.3263)}      & \color{black} 5.27         & \color{black}1288.41  & \color{black}1440.3    \\
  \hline\hline
  \end{tabular}
  \begin{tablenotes}
  \footnotesize
\item[*] IGBT-R provides reactive power in \textcircled{2}-\textcircled{5} and operates at PF = 1 in \textcircled{6}-\textcircled{9}.
\end{tablenotes}
\end{threeparttable}
  \vspace{0pt}
\end{table}

\subsubsection{\color{black}Comparisons with other planning methods and existing configurations}
\label{sec:M1/M2}

{\color{black}
To verify the effectiveness of incorporating active-reactive power coordination and the coupling between hydrogen production and the electrical network into the planning process, and considering the absence of such integrated methods in existing literature, the proposed method is compared against the following intuitive schemes:
\begin{itemize}
  \item \textbf{M1}: Focusing only on hydrogen plant performance while neglecting the reactive load of rectifiers.  This approach is commonly adopted in existing hydrogen plant planning studies \cite{li2019optimal, li2023exploration, varela2021modeling}.
  \item\textbf{M2}: Considers the reactive load but does not account for reactive power exchange with the network; instead, the SVG capacity is sized based on the peak reactive demand of the hydrogen plant.s
  \item \textbf{M3}/\textbf{M4}: The number of ELZs is fixed to match an actual project \cite{huadian2024} ($\sigma^{\text{ELZ}}=3$), while other planning settings for rectifier configuration are kept the same as in M1/M2, respectively.
\end{itemize}
}

{\color{black}
The planning results are summarized in Table \ref{tab:M1}. Compared with M1-M4, the proposed method increases total revenue by up to 13.78\%. In all M1-M4 cases, SVG capacity is not jointly optimized, and TRs are selected exclusively due to their higher efficiency. However, in M1 and M3, no SVG is deployed to compensate for the reactive demand of TRs. To maintain reactive power balance, the system must rely on alternative sources (e.g., WT, PV, BES) to supply reactive power, leading to excessive var flow in the network. This not only reduces the available capacity for active power delivery, but also exacerbates voltage security issues and increases network losses, ultimately lowering RE utilization and system-level economic performance.
In contrast, M2 and M4 neglect the reactive power flexibility of the hydrogen plant, resulting in oversized SVG capacity and poor cost-effectiveness. These results collectively demonstrate that coordinated planning of active and reactive power flows is critical to improving the techno-economic performance of ReP2H systems and should therefore be considered essential.
}


\begin{table}[t]\scriptsize \color{black}
  \renewcommand{\arraystretch}{1.45}
  \caption{Comparison Between the Proposed Method and Other Planning Methods and Configurations in the Small-Scale System}\vspace{-0pt}
  \label{tab:M1}
  \centering
  \begin{tabular}{c@{\hspace{5pt}}c@{\hspace{5pt}}c@{\hspace{5pt}}c@{\hspace{5pt}}c@{\hspace{5pt}}c}
  \hline\hline
   Method       & Planning Results           & \tabincell{c}{Total Revenue\\ (10$^7$ CNY/yr)}   & \tabincell{c}{Improvement of the\\ Proposed Method}    \\
  \hline
  \textbf{M1}    &  4 TRs + 0 MVar SVG        &  1.8267     & \textbf{+13.78\%}     \\
  \textbf{M2}    &  4 TRs + 14.36 MVar SVG    &  2.0468    & \textbf{+1.55\%}     \\
  \textbf{M3}    &  3 TRs + 0 MVar SVG        &  1.8471         & \textbf{+12.53\%}    \\
  \textbf{M4}    &  3 TRs + 10.77 MVar SVG      &  2.0375     &  \textbf{+ 2.01\%}      \\
  \hline
  \textbf{Proposed}       &  Configuration \textcircled{2}   & 2.0785     &     \\
  \hline\hline
  \end{tabular}
\end{table}

\subsection{Operational Analysis of Different Rectifier Configurations 
}
\label{sec:Comparison}

\subsubsection{\color{black}Scheduling strategy of ELZs}

\begin{figure}[t]   \centering
  \includegraphics[width=3.4in]{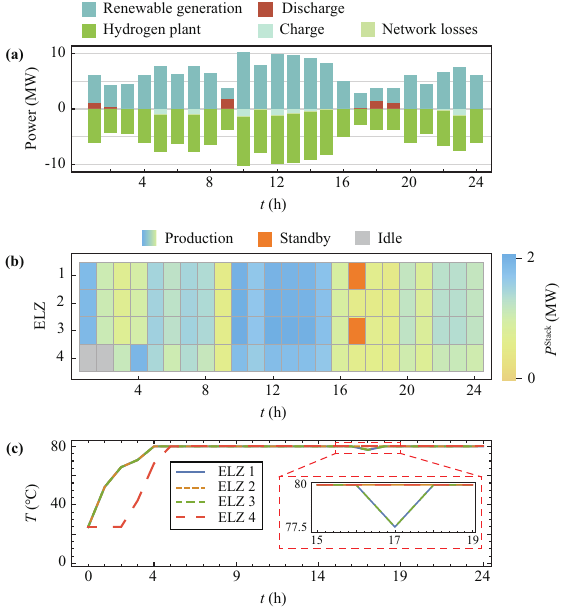}\vspace{-0pt}
  \caption{{\color{black}Scheduling strategy of the hydrogen plant under the optimal configuration \textcircled{2} (3 TRs and 1 IGBT-R) in a specific scenario. (a) Power balance. (b) State switches and electrolytic loads. (c) Stack temperature control.}}
  \label{fig:scheduling}\vspace{-0pt}
\end{figure}

{\color{black}
We analyze power balance, ELZ state transitions, load allocation, and temperature control to demonstrate the core logic and strategy of active power scheduling for hydrogen production. The operational results for a specific scenario using the optimal configuration \textcircled{2} are illustrated in Fig. \ref{fig:scheduling}.

As shown in Fig. \ref{fig:scheduling}(a), hydrogen production closely follows the fluctuations in renewable generation, with the BES providing load balancing and enabling continuous ELZ operation. The hydrogen plant achieves wide-range power regulation through dynamic switching among \textit{production}, \textit{standby}, and \textit{idle} states across the ELZs. Notably, the 4th ELZ operates at a lower load due to the lower efficiency of its rectifier, while the other three ELZs adhere to the equimarginal principle \cite{zeng2024scheduling,li2024two} to maximize hydrogen yield.
In addition, thermal management plays a vital role. At $t = \{3, 4\}$ h, the 4th ELZ operates at the highest power to use electrolytic heat for P2H efficiency improvement. For details on network-side operations, interested readers are referred to \cite{zeng2024scheduling}, and the reactive power analysis will be discussed later.

}

\subsubsection{\color{black}Operational comparison}

\begin{figure}[t]
  \centering
  \includegraphics[width=3.35in]{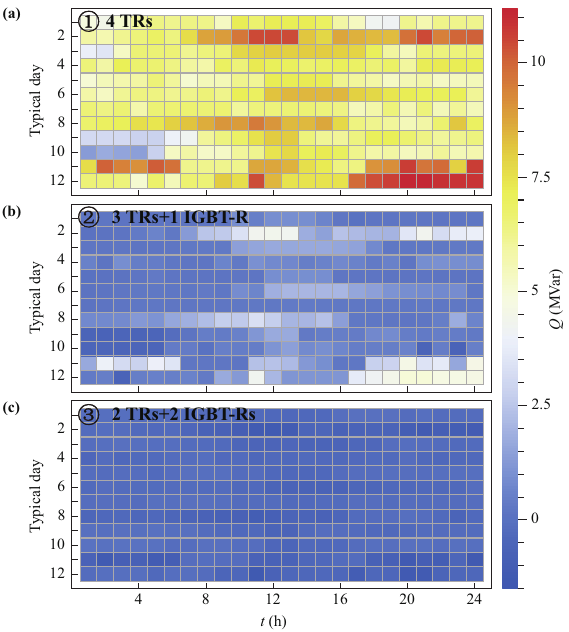}\vspace{-0pt}
  \caption{Reactive power demand of the hydrogen plant with different rectifier configurations under 12 scenarios in the first representative year. (a) \textcircled{1} 4 TRs. (b) \textcircled{2} 3 TRs + 1 IGBT-R. (c) \textcircled{3} 2 TRs + 2 IGBT-Rs.}
  \label{fig:reactivepower}\vspace{-0pt}
\end{figure}

\begin{figure}[t]   \centering
  \includegraphics[width=3.35in]{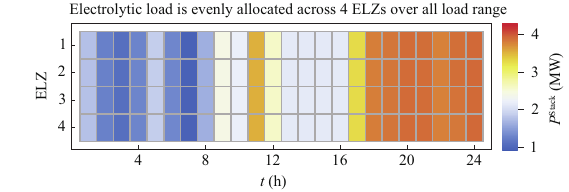}\vspace{-0pt}
  \caption{Electrolytic load of 4 ELZs in the typical renewable generation scenario under configuration \textcircled{1} with 4 TRs.}
  \label{fig:PQ1}\vspace{6pt}
  \includegraphics[width=3.35in]{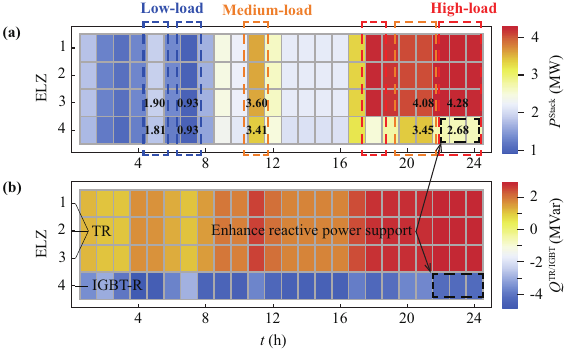}\vspace{-0pt}
  \caption{(a) DC-side electrolytic load and (b) AC-side reactive power of 4 ELZs in the typical scenario under the optimal configuration \textcircled{2} with 3 TRs and 1 IGBT-R.}
  \label{fig:PQ2}\vspace{6pt}
  \includegraphics[width=3.35in]{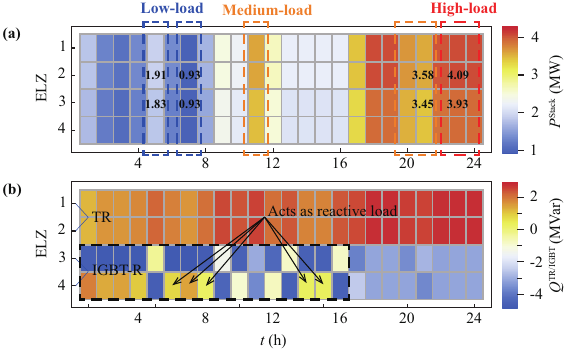}\vspace{-0pt}
  \caption{(a) DC-side electrolytic load and (b) AC-side reactive power of 4 ELZs in the typical scenario under configuration \textcircled{3} with 2 TRs and 2 IGBT-Rs.}
  \label{fig:PQ3}\vspace{-0pt}
\end{figure}

To analyze the operational differences of the ReP2H system under various rectifier configurations, we focus on configurations \textcircled{1}-\textcircled{3}. Configurations \textcircled{4} and \textcircled{5} are similar to \textcircled{3}, whereas \textcircled{6}-\textcircled{9} lack economic advantages and are thus left out here. The total reactive power demand of the hydrogen plant under different configurations and typical days (for a single hydrogen price scenario) is shown in Fig. \ref{fig:reactivepower}. The electrolytic load and reactive power of the four ELZs in a typical renewable generation scenario under configurations \textcircled{1}-\textcircled{3} are illustrated in Figs. \ref{fig:PQ1} to \ref{fig:PQ3}, respectively.

Fig. \ref{fig:reactivepower} shows that a hydrogen plant powered by four TRs exhibits significant reactive power demand (over 6 MVar) across most scenarios. In contrast, replacing one TR with an IGBT-R can reduce the reactive power demand to below 1.5 MVar, achieving local balance and reducing the reactive power flow and network losses.
Configuration \textcircled{3} with 2 TRs and 2 IGBT-Rs provides reactive power support of approximately 0--0.6 MVar to the network. Clearly, this configuration has excessive var resources, making it less cost-effective.

As shown in Fig. \ref{fig:PQ1}, when the rectifiers in all ELZs are uniform and var resources are adequate, the electrolytic load is evenly distributed across the ELZs to maximize the hydrogen yield, aligning with the equimarginal principle and resulting in \cite{zeng2024scheduling}. In contrast, mixed rectifier configurations \textcircled{2} and \textcircled{3}, as shown in Figs. \ref{fig:PQ2} and \ref{fig:PQ3}, result in ELZs powered by TRs bearing a slightly higher load than those powered by IGBT-Rs in \textcircled{2} and \textcircled{3}.

Fig. \ref{fig:PQ2}(a) shows that the load difference between the TR- and IGBT-powered ELZs is the smallest at low-load intervals, increases at medium-load intervals, and peaks at high-load intervals.
Combining Fig. \ref{fig:PQ2}(b), we see that reducing the load of IGBT-powered ELZs releases their reactive support capability (see Fig. \ref{fig:region}(b)), facilitating local reactive power balance at the hydrogen plant at high loads, minimizing network losses, and ensuring voltage stability. However, the concentrating load in a few ELZs lowers the overall P2H energy conversion efficiency, as dictated by the electrochemical characteristics (Section \ref{sec:same}) and the equimarginal principle \cite{li2022coordinated,zeng2024scheduling}. Thus, during low- and medium-load intervals, when the reactive power demand is lower, a more even load distribution improves the hydrogen yield. Furthermore, Fig. \ref{fig:PQ3}(b) shows that at low- and medium-load intervals, the IGBT-R even acts as a reactive load, confirming a surplus in var resources in configuration \textcircled{3} (2 TRs, 2 IGBT-Rs).


\subsection{Robust Portfolio Results of the Small-scale System for\\$\beta>0$ }

\begin{figure}[t]
	\centering
	\includegraphics[width=3.4in]{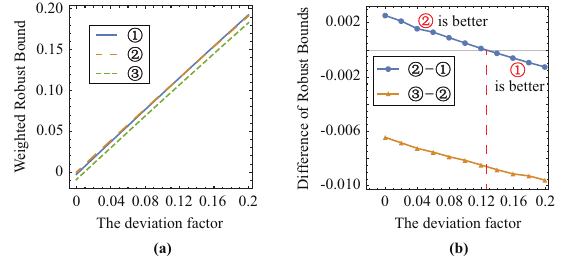}\vspace{-0pt}
	\caption{Weighted robust uncertainty bounds of different $\beta$ under configurations \textcircled{1} 4 TRs, \textcircled{2} 3 TRs + 1 IGBT-R and \textcircled{3} 2 TRs + 2 IGBT-Rs, and the differences among them. (a) Robust bound. (b) Difference.}
	\label{fig:IGDT}\vspace{-0pt}
\end{figure}

With four ELZs and a revenue deviation factor $\beta$ ranging from 0 to 0.2, Fig. \ref{fig:IGDT}(a) presents the maximum weighted robust uncertainty bounds for different configurations, whereas Fig. \ref{fig:IGDT}(b) illustrates the differences among them.

As depicted in Fig. \ref{fig:IGDT}(a), decreasing $\beta$ increases both the robust bound and the acceptable risk. In managing renewable generation uncertainty, configurations \textcircled{1} and \textcircled{2} perform similarly and outperform \textcircled{3}. Fig. \ref{fig:IGDT}(b) further shows that as $\beta$ increases, the advantage of configuration \textcircled{2} over \textcircled{1} diminishes. This is because the reduced renewable output relaxes the network power flow constraints, leading to a reduction in the reactive power demand. When $\beta > 0.12$, configuration \textcircled{1} outperforms \textcircled{2}. This demonstrates that under high uncertainty risk, the IGDT model provides better results than approaches that do not consider intraday uncertainties. With respect to the rectifier configuration, a higher level of renewable output and a higher proportion of stressed operations increase the demand for IGBT-Rs, whereas a lower output level favors TRs.

{\color{black}
	\subsection{Sensitivity Analysis}
	\label{sec:sensitivity}
	
	\subsubsection{BES capacity}

	As the BES participates in both active and reactive power scheduling, it is necessary to examine whether its capacity influences system operation and, consequently, the rectifier portfolio. To this end, a sensitivity analysis is conducted to assess how BES capacity affects the optimal rectifier configuration. The corresponding system revenue is presented in Fig. \ref{fig:BESsensi}(a).
	The results show that system revenue increases with BES capacity. However, the relative performance of rectifier configurations remains consistent: $\textcircled{2} > \textcircled{1} > \textcircled{3}$, with the MR configuration being the most favorable, in line with the baseline result. This indicates that within a reasonable capacity range, the BES primarily influences active power balancing, while its impact on reactive power flow or the optimal rectifier configuration is negligible due to location differences.
	
	\begin{figure}[t]
		\centering
		\includegraphics[width=3.35in]{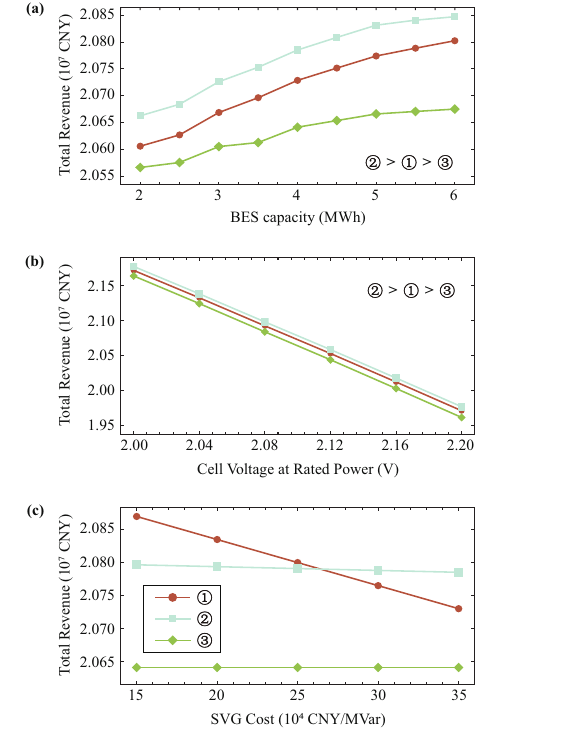}\vspace{-0pt}
		\caption{{\color{black}Impact of (a) BES Capacity, (b) Stack U-I Characteristics, and (c) SVG cost on system revenue and optimal rectifier configuration.}}
		\label{fig:BESsensi}\vspace{-0pt}
	\end{figure}
	
	\subsubsection{Stack U-I characteristics}
	
	Due to manufacturing variability and stack degradation, the performance of electrolysis stacks can vary, primarily reflected in their U-I characteristics. To evaluate the impact of such deviations, we examine rectifier planning under different stack performances. At rated power, cell voltage ranges from 2.0 to 2.2 V (with 2.1 V as the baseline), where lower voltage indicates higher efficiency.
	Fig. \ref{fig:BESsensi}(b) shows system revenue under various stack performance levels. It is evident that configuration \textcircled{2} consistently yields the highest revenue. Thus, although stack performance affects overall efficiency, its impact on the relative performance of rectifier configurations is uniform and does not alter the optimal configuration ratio.

	\subsubsection{SVG cost}
	
	To further identify factors influencing rectifier selection, we analyze the effect of SVG cost. Fig. \ref{fig:BESsensi}(c) presents the variation in system revenue with respect to SVG cost. As SVG becomes cheaper, the revenue associated with the 4-TR configuration increases noticeably, while configurations involving IGBT-Rs are only marginally affected. When SVG cost drops to approximately 250 CNY/kVar, the 4-TR setup becomes the most cost-effective option.
	
	In summary, small variations in most system parameters do not affect the optimal MR configuration, as the self-balancing capability of reactive power within the hydrogen plant ensures an already efficient state. However, SVG cost plays a significant role in economic performance while maintaining a reactive power balance, making it a critical factor in rectifier configuration planning.
	
}

\subsection{Application in a Large-Scale System}
\label{sec:caselarge}

Finally, to verify scalability and practical applicability, we apply the proposed model to the large-scale system briefly described in Section \ref{sec:largecase}. A 4-in-1 ELZ configuration is adopted, where four ELZs share a set of BoPs (e.g., lye-gas separators, lye circulation loops, and purification units) to reduce the footprint and cost \cite{liang2024large, qiu2025dynamic}. Thus, the production-standby-idle states and temperature dynamics of each 4-in-1 unit are coupled, while the active and reactive powers remain independently scheduled.

\begin{figure}[t]
  \centering
  \includegraphics[width=3.4in]{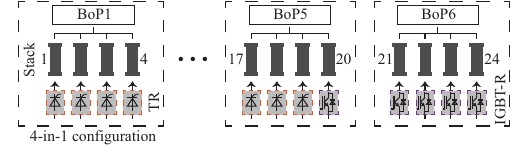}\vspace{-0pt}
  \caption{Detailed configurations of the 24 ELZs in the large-scale system.}   \label{fig:4in1}\vspace{-0pt}
\end{figure}

\begin{figure}[t]
  \centering
  \includegraphics[width=3.45in]{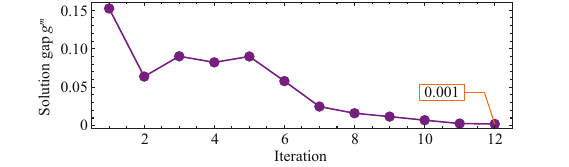}\vspace{-0pt}
  \caption{{\color{black}Convergence of the PH algorithm for solving the SP model.}}   \label{fig:convergence}\vspace{-0pt}
\end{figure}

For a deviation factor $\beta \in[0,0.06]$, the optimal rectifier portfolio consists of 24 ELZs powered by 19 TRs and 5 IGBT-Rs, shown in Fig. \ref{fig:4in1}. Notably, the positioning of TRs and IGBT-Rs within the 4-in-1 ELZs, i.e., whether uniform or mixed, is not optimized here and will be explored in future work. {\color{black}In terms of computational performance, the simulation of a large-scale case is performed on a desktop with an \textit{Intel Core i5-12400@2.5 GHz} CPU and 16 GB of RAM, and solved by \textit{Gurobi 11.0}. The convergence of the PH algorithm applied to the SP model is shown in Fig. \ref{fig:convergence}, with a computation time of approximately 24 minutes. This indicates that the large-scale planning problem in this study can be efficiently solved after decomposition, demonstrating the scalability.}

Table \ref{tab:compasison2} compares the total revenue for the optimal rectifier configuration, uniform TR/IGBT-R configurations, ratios used in existing projects, {\color{black}and configurations planned by M1 and M2 (introduced in Section \ref{sec:M1/M2})} for $\beta=0$. Compared with the uniform configurations, the optimal configuration increases total revenue by 0.89\% to 2.56\%, as well as 0.79\%, 1.37\%, compared with the 1:1 (32 TRs + 32 IGBT-Rs in the Songyuan Project \cite{songyuan2023}) and 1:2 (12 TRs + 24 IGBT-Rs in the Da'an Project \cite{daan2024}) configurations.
{\color{black}
	Most importantly, the proposed method increases total revenue by at most 10.70\% compared to the methods M1 and M2 used in the existing hydrogen planning studies without considering the reactive power and the network. We also find that}
the suboptimal portfolio (18 TRs, 6 IGBT-Rs) achieves nearly identical revenue to the optimal portfolio.

\begin{table}[t]\scriptsize
  \renewcommand{\arraystretch}{1.45}
  \caption{Comparison of Different Rectifier and Var Compensation Configurations in the Large-Scale System for $\beta=0$}\vspace{-0pt}
  \label{tab:compasison2}
  \centering
  \begin{tabular}{c@{\hspace{5pt}}c@{\hspace{5pt}}c@{\hspace{5pt}}c}
  \hline\hline
   Configuration                     &\tabincell{c}{ Total Revenue\\(10$^8$ CNY/yr)} & \tabincell{c}{SVG Capacity\\(MVar)}    & \tabincell{c}{Comparision} \\
  \hline
  24 TRs (uniform)                  & 1.2248                  & 44     &  \textbf{+0.89\%}            \\
  18 TRs + 6 IGBT-Rs (3:1)          & 1.2357                  & 0     &  \textbf{$\approx$0}            \\
  12 TRs + 12 IGBT-Rs (1:1)         & 1.2261                  & 0     &  \textbf{+0.79\%}            \\
  8 TRs + 16 IGBT-Rs (1:2)          & 1.2191                  & 0     &  \textbf{+1.37\%}               \\
  24 IGBT-Rs (uniform)              & 1.2049                 & 0      &  \textbf{+2.56\%}             \\
  24 \color{black}  TRs (M1)         & \color{black}  1.1164                  &\color{black}   0      &\color{black}   \textbf{+10.70\%}            \\
  24 \color{black}  TRs (M2)             &\color{black}  1.2113                  &\color{black}  86.16     &\color{black}   \textbf{+2.02\%}            \\
  \hline
  19 TRs + 5 IGBT-Rs (optimal)      & \textbf{1.2358}         & 0           &               \\
  \hline\hline
  \end{tabular}\vspace{-0pt}
\end{table}

\begin{figure}[t]   \centering
  \includegraphics[width=3.4in]{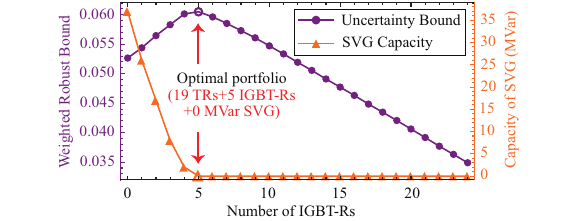}\vspace{-0pt}
  \caption{Weighted robust uncertainty bounds and optimal SVG capacities of different rectifier configurations in the large-scale system for $\beta=0.06$.}
  \label{fig:revenue2}\vspace{-0pt}
\end{figure}

Fig. \ref{fig:revenue2} displays the robust results for $\beta=0.06$ under different rectifier configurations. The optimal SVG capacity remains 0 MVar, which is consistent with $\beta=0$, but the suboptimal configuration shifts to 20 TRs and 4 IGBT-Rs. This further confirms that IGBT-R performs better in stressed scenarios, with the full load hours of renewable energy or ELZs being key factors influencing the TR/IGBT-R ratio.

From Tables \ref{tab:planning} and \ref{tab:compasison2} and Fig. \ref{fig:revenue2}, we can summarize that the optimality criterion for MR configurations can be understood as balancing the reactive power in the hydrogen plant, optimizing the investment in SVG to approximately 0, or fully utilizing the reactive power support capability of IGBT-Rs. Excessive investment in IGBT-Rs leads to idle var resources and reduced profitability. Moreover, for low renewable uncertainty, a portfolio of 18 TRs and 6 IGBT-Rs remains near-optimal. Thus, we conclude that the ideal TR-to-IGBT-R ratio for electrolysis rectifiers is approximately 3:1.

\section{Conclusions}
\label{sec:conclusion}

This paper examines the complementarity of TR and IGBT-R in powering electrolyzers within ReP2H systems, explores the feasibility of mixed rectifier configurations, and integrates them into an investment planning framework. An optimal portfolio planning model and its solution approach are proposed. Simulations based on real-world systems in Inner Mongolia, China,
yield the following key findings:

1) In the optimal portfolio, the reactive power support capability of the IGBT-Rs must be fully utilized {\color{black}to achieve self-balancing of reactive power within the hydrogen plant}, while the need for var compensation approaches zero.

2) The ideal TR-to-IGBT-R ratio for electrolysis rectifiers in ReP2H systems is approximately 3:1 {\color{black}given the current costs of rectifiers and SVG}.

3) The renewable energy profile influences the rectifier portfolio, as it determines the network power flow stress. A higher proportion of stressed operations increases the demand for IGBT-Rs, whereas a lower proportion favors TRs.

Future research will focus on the transient behavior and performance of different electrolysis rectifiers, as well as their interactive mechanisms with the electrical network, to enhance the strength and operational stability of ReP2H systems {\color{black}both in steady-state operation and ride-through events}, particularly under weak or off-grid conditions.


\begin{thebibliography}{10}
\providecommand{\url}[1]{#1}
\csname url@samestyle\endcsname
\providecommand{\newblock}{\relax}
\providecommand{\bibinfo}[2]{#2}
\providecommand{\BIBentrySTDinterwordspacing}{\spaceskip=0pt\relax}
\providecommand{\BIBentryALTinterwordstretchfactor}{4}
\providecommand{\BIBentryALTinterwordspacing}{\spaceskip=\fontdimen2\font plus
\BIBentryALTinterwordstretchfactor\fontdimen3\font minus
  \fontdimen4\font\relax}
\providecommand{\BIBforeignlanguage}[2]{{%
\expandafter\ifx\csname l@#1\endcsname\relax
\typeout{** WARNING: IEEEtran.bst: No hyphenation pattern has been}%
\typeout{** loaded for the language `#1'. Using the pattern for}%
\typeout{** the default language instead.}%
\else
\language=\csname l@#1\endcsname
\fi
#2}}
\providecommand{\BIBdecl}{\relax}
\BIBdecl

\bibitem{van2020hydrogen}
S.~Van~Renssen, ``The hydrogen solution?'' \emph{Nature Clim. Chang.}, vol.~10,
  no.~9, pp. 799--801, 2020.

\bibitem{guo2023deploying}
Y.~Guo, L.~Peng, J.~Tian, and D.~L. Mauzerall, ``Deploying green hydrogen to
  decarbonize {C}hina's coal chemical sector,'' \emph{Nature Commun.}, vol.~14,
  no.~1, p. 8104, 2023.

\bibitem{odenweller2022probabilistic}
A.~Odenweller, F.~Ueckerdt, G.~F. Nemet, M.~Jensterle, and G.~Luderer,
  ``Probabilistic feasibility space of scaling up green hydrogen supply,''
  \emph{Nature Energy}, vol.~7, no.~9, pp. 854--865, 2022.

\bibitem{2022}
\BIBentryALTinterwordspacing
{The Energy Bureau of Inner Mongolia Autonomous Region}, ``Notice of {Inner
  Mongolia Autonomous Region Energy Bureau} on carrying out the 2022 wind-solar
  hydrogen production integration demonstration.'' [Online]. Available:
  \url{http://dbnyb.com/07/taiyangnen/2021/0827/51472.html}
\BIBentrySTDinterwordspacing

\bibitem{qiu2023extend}
Y.~Qiu, B.~Zhou, T.~Zang, Y.~Zhou, S.~Chen, R.~Qi, J.~Li, and J.~Lin,
  ``Extended load flexibility of utility-scale {P2H} plants: Optimal production
  scheduling considering dynamic thermal and {HTO} impurity effects,''
  \emph{Renewable Energy}, vol. 217, p. 119198, Nov. 2023.

\bibitem{zeng2024scheduling}
Y.~Zeng, Y.~Qiu, J.~Zhu, S.~Chen, B.~Zhou, J.~Li, B.~Yang, and J.~Lin,
  ``Scheduling multiple industrial electrolyzers in renewable {P2H} systems: A
  coordinated active-reactive power management method,'' \emph{IEEE Trans.
  Sustain. Energy}, vol.~16, no.~1, pp. 201--215, Jan. 2025.

\bibitem{koponen2021comparison}
J.~Koponen, A.~Poluektov, V.~Ruuskanen, A.~Kosonen, M.~Niemel{\"a}, and
  J.~Ahola, ``Comparison of thyristor and insulated-gate bipolar
  transistor-based power supply topologies in industrial water electrolysis
  applications,'' \emph{J. Power Sources}, vol. 491, p. 229443, Apr. 2021.

\bibitem{songyuan2023}
\BIBentryALTinterwordspacing
{China Energy Engineering}, ``{C}hina {E}nergy {E}ngineering's {S}ongyuan
  hydrogen industry park,'' 2023. [Online]. Available:
  \url{http://www.ceehe.ceec.net.cn/art/2023/9/26/art\_59190\_802.html}
\BIBentrySTDinterwordspacing

\bibitem{daan2024}
\BIBentryALTinterwordspacing
{Da'an Municipal People's Government}, ``Da'an integrated green hydrogen and
  ammonia demonstration project: Building momentum for a green future,'' 2024.
  [Online]. Available:
  \url{http://daan.jlbc.gov.cn/daxw/ztbd/202410/t20241007\_1000080.html}
\BIBentrySTDinterwordspacing

{\color{black}\bibitem{wu2024coordination}
Z.~Wu, ``Coordination control of smart electrolyzer and grid based on
  current-source rectifier under {LVRT},'' in \emph{2024 3rd Asia Power and
  Electrical Technology Conference}, 2024, pp. 498--502.

\bibitem{meng2022novel}
X.~Meng, M.~Chen, M.~He, X.~Wang, and J.~Liu, ``A novel high power hybrid
  rectifier with low cost and high grid current quality for improved efficiency
  of electrolytic hydrogen production,'' \emph{IEEE Trans. Power Electron.},
  vol.~37, no.~4, pp. 3763--3768, Apr. 2022.}

\bibitem{de2023hydrogen}
A.~M. De~Corato, M.~Ghazavi~Dozein, S.~Riaz, and P.~Mancarella, ``Hydrogen
  electrolyzer load modelling for steady-state power system studies,''
  \emph{IEEE Trans. Power Deliv.}, vol.~38, no.~6, pp. 4312--4323, Dec. 2023.

\bibitem{zhang2022research}
M.~Zhang, R.~Liu, K.~Wang, K.~Sun, Q.~Guo, and Y.~Li, ``Research on low voltage
  ride through and reactive power support of hydrogen production power
  supply,'' in \emph{2022 4th Int. Conf. Smart Power \& Internet Energy Syst.},
  2022, pp. 438--442.

\bibitem{tavakoli2023gridforming}
S.~D. Tavakoli, M.~G. Dozein, V.~A. Lacerda, M.~C. Mañe, E.~Prieto-Araujo,
  P.~Mancarella, and O.~Gomis-Bellmunt, ``Grid-forming services from hydrogen
  electrolyzers,'' \emph{IEEE Trans. Sustain. Energy}, vol.~14, no.~4, pp.
  2205--2219, Oct. 2023.

\bibitem{ruuskanen2020power}
V.~Ruuskanen, J.~Koponen, A.~Kosonen, M.~Niemel{\"a}, J.~Ahola, and
  A.~H{\"a}m{\"a}l{\"a}inen, ``Power quality and reactive power of water
  electrolyzers supplied with thyristor converters,'' \emph{J. Power Sources},
  vol. 459, p. 228075, May 2020.

\bibitem{li2024two}
J.~Li, B.~Yang, J.~Lin, F.~Liu, Y.~Qiu, Y.~Xu, R.~Qi, and Y.~Song, ``Two-layer
  energy management strategy for grid-integrated multi-stack power-to-hydrogen
  station,'' \emph{Appl. Energy}, vol. 367, p. 123413, Aug. 2024.

\bibitem{gao2024advanced}
Y.~Gao, X.~Wang, and X.~Meng, ``Advanced rectifier technologies for
  electrolysis-based hydrogen production: A comparative study and real-world
  applications,'' \emph{Energies}, vol.~18, no.~1, p.~48, Dec. 2024.

\bibitem{intellipower2022}
\BIBentryALTinterwordspacing
Intellipower, ``Harmonic analysis and solutions for large-scale water
  electrolysis hydrogen production,'' 2022. [Online]. Available:
  \url{https://www.intellipower.cn/industry-735}
\BIBentrySTDinterwordspacing

\bibitem{li2019optimal}
J.~Li, J.~Lin, H.~Zhang, Y.~Song, G.~Chen, L.~Ding, and D.~Liang, ``Optimal
  investment of electrolyzers and seasonal storages in hydrogen supply chains
  incorporated with renewable electric networks,'' \emph{IEEE Trans. Sustain.
  Energy}, vol.~11, no.~3, pp. 1773--1784, Jul. 2020.

\bibitem{pan2020optimal}
G.~Pan, W.~Gu, Y.~Lu, H.~Qiu, S.~Lu, and S.~Yao, ``Optimal planning for
  electricity-hydrogen integrated energy system considering power to hydrogen
  and heat and seasonal storage,'' \emph{IEEE Trans. Sustain. Energy}, vol.~11,
  no.~4, pp. 2662--2676, Oct. 2020.

\bibitem{zhu2024full}
J.~Zhu, Y.~Qiu, Y.~Zeng, Y.~Zhou, S.~Chen, T.~Zang, B.~Zhou, Z.~Yu, and J.~Lin,
  ``Exploring the optimal size of grid-forming energy storage in an off-grid
  renewable {P2H} system under multi-timescale energy management,'' \emph{arXiv
  preprint arXiv:2409.05086}, 2024.

\bibitem{li2023exploration}
Y.~Li, X.~Deng, T.~Zhang, S.~Liu, L.~Song, F.~Yang, M.~Ouyang, and X.~Shen,
  ``Exploration of the configuration and operation rule of the
  multi-electrolyzers hybrid system of large-scale alkaline water hydrogen
  production system,'' \emph{Appl. Energy}, vol. 331, p. 120413, Feb. 2023.

\bibitem{ibanez2023off}
A.~Ib{\'a}{\~n}ez-Rioja, L.~J{\"a}rvinen, P.~Puranen, A.~Kosonen, V.~Ruuskanen,
  K.~Hynynen, J.~Ahola, and P.~Kauranen, ``Off-grid solar {PV}--wind
  power--battery--water electrolyzer plant: Simultaneous optimization of
  component capacities and system control,'' \emph{Appl. Energy}, vol. 345, p.
  121277, Sep. 2023.

\bibitem{zheng2023model}
Y.~Zheng, S.~You, C.~Huang, and X.~Jin, ``Model-based economic analysis of
  off-grid wind/hydrogen systems,'' \emph{Renewable Sustain. Energy Rev.}, vol.
  187, p. 113763, Nov. 2023.

\bibitem{escn2024}
\BIBentryALTinterwordspacing
{ESCN}, ``Intelli {L}i {H}uan: Multiple hydrogen production rectifier
  solutions, respecting the first principles of technology,'' 2024. [Online].
  Available: \url{https://www.escn.com.cn/20240827/
  d9f3804b25104b8ba5101e0e410c6ac0/c.html}
\BIBentrySTDinterwordspacing

{\color{black}\bibitem{wang2024optimization}
X.~Wang, X.~Meng, G.~Nie, B.~Li, H.~Yang, and M.~He, ``Optimization of hydrogen
  production in multi-electrolyzer systems: A novel control strategy for
  enhanced renewable energy utilization and electrolyzer lifespan,''
  \emph{Appl. Energy}, vol. 376, p. 124299, Dec. 2024.}

\bibitem{liang2024large}
T.~Liang, M.~Chen, J.~Tan, Y.~Jing, L.~Lv, and W.~Yang, ``Large-scale off-grid
  wind power hydrogen production multi-tank combination operation law and
  scheduling strategy taking into account alkaline electrolyzer
  characteristics,'' \emph{Renewable Energy}, vol. 232, p. 121122, Oct. 2024.

{\color{black}\bibitem{firdous2025short}
A.~Firdous, C.~P. Barala, P.~Mathuria, and R.~Bhakar, ``Short-term operation
  flexibility in modular power to hydrogen based ammonia industries,''
  \emph{IEEE Trans. Sustain. Energy}, vol.~16, no.~1, pp. 601--612, Jan. 2025.

\bibitem{Khajeh2024Optimized}
H.~Khajeh, S.~Seyyedeh-Barhagh, and H.~Laaksonen, ``Optimized operation of
  hybrid wind-hydrogen system to provide flexibility for transmission system
  needs,'' \emph{IEEE Trans. Sustain. Energy}, pp. 1--13, 2024, Early Access.

\bibitem{Han2025Robust}
P.~Han, X.~Xu, Z.~Yan, M.~Shahidehpour, Z.~Tan, H.~Wang, and G.~Li, ``Robust
  scheduling of integrated electricity-heat-hydrogen system considering
  bidirectional heat exchange between alkaline electrolyzers and district
  heating networks,'' \emph{J. Mod. Power Syst. Clean Energy}, pp. 1--13, 2025, Early Access.}

\bibitem{rockafellar1991scenarios}
R.~T. Rockafellar and R.~J.-B. Wets, ``Scenarios and policy aggregation in
  optimization under uncertainty,'' \emph{Math. Oper. Res.}, vol.~16, no.~1,
  pp. 119--147, Feb. 1991.

\bibitem{neimenggu2024}
\BIBentryALTinterwordspacing
{The Energy Bureau of Inner Mongolia Autonomous Region}, ``Implementation
  guidelines for wind and solar power to ammonia and methanol production
  projects in {Inner Mongolia Autonomous Region},'' 2024. [Online]. Available:
  \url{https://nyj.nmg.gov.cn/zwgk/zfxxgkzl/fdzdgknr/tzgg\_16482/gg\_16484/202409
  /t20240903\_2568246.html}
\BIBentrySTDinterwordspacing

\bibitem{ulleberg2003modeling}
{\O}.~Ulleberg, ``Modeling of advanced alkaline electrolyzers: a system
  simulation approach,'' \emph{Int. J. Hydrogen Energy}, vol.~28, no.~1, pp.
  21--33, Jan. 2003.

\bibitem{varela2021modeling}
C.~Varela, M.~Mostafa, and E.~Zondervan, ``Modeling alkaline water electrolysis
  for power-to-x applications: A scheduling approach,'' \emph{Int. J. Hydrogen
  Energy}, vol.~46, no.~14, pp. 9303--9313, Feb. 2021.

\bibitem{biaozhun}
{China Electricity Council}, ``Technical specification for rectifier power
  supply in water electrolysis hydrogen production,'' Standard No. T/CES
  226-2023, 2023, in Chinese.

{\color{black}\bibitem{jiang2025two}
X.~Jiang and Q.~Xu, ``Two-layer iterative sizing for logistics centers equipped
  with fast charging stations,'' \emph{IEEE Trans. Smart Grid}, 2025, Early Access.}

\bibitem{li2022coordinated}
J.~Li, J.~Lin, Y.~Song, J.~Xiao, F.~Liu, Y.~Zhao, and S.~Zhan, ``Coordinated
  planning of {HVDC}s and power-to-hydrogen supply chains for interregional
  renewable energy utilization,'' \emph{IEEE Trans. Sustain. Energy}, vol.~13,
  no.~4, pp. 1913--1929, Oct. 2022.

\bibitem{yang2024stochastic}
L.~Yang, H.~Li, H.~Zhang, Q.~Wu, and X.~Cao, ``Stochastic-distributionally
  robust frequency-constrained optimal planning for an isolated microgrid,''
  \emph{IEEE Trans. Sustain. Energy}, vol.~15, no.~4, pp. 2155--2169, Oct.
  2024.

\bibitem{yu2024optimal}
Z.~Yu, J.~Lin, F.~Liu, J.~Li, Y.~Zhao, Y.~Song, Y.~Song, and X.~Zhang,
  ``Optimal sizing and pricing of grid-connected renewable power to ammonia
  systems considering the limited flexibility of ammonia synthesis,''
  \emph{IEEE Trans. Power Syst.}, vol.~39, no.~2, pp. 3631--3648, Mar. 2024.

\bibitem{liu2024system}
Y.~Liu, Y.~Chen, H.~Xin, J.~Tu, L.~Zhang, M.~Song, and J.~Zhu, ``System
  strength constrained grid-forming energy storage planning in renewable power
  systems,'' \emph{IEEE Trans. Sustain. Energy}, pp. 1--13, 2024.

\bibitem{farivar2013branch}
M.~Farivar and S.~H. Low, ``Branch flow model: Relaxations and
  convexification---{P}art {I},'' \emph{IEEE Trans. Power Syst.}, vol.~28,
  no.~3, pp. 2554--2564, Aug. 2013.

\bibitem{ghazavi2021simultaneous}
M.~Ghazavi~Dozein, O.~Gomis-Bellmunt, and P.~Mancarella, ``Simultaneous
  provision of dynamic active and reactive power response from utility-scale
  battery energy storage systems in weak grids,'' \emph{IEEE Trans. Power
  Syst.}, vol.~36, no.~6, pp. 5548--5557, Nov. 2021.

\bibitem{nasr2020assessing}
M.-A. Nasr, E.~Nasr-Azadani, H.~Nafisi, S.~H. Hosseinian, and P.~Siano,
  ``Assessing the effectiveness of weighted information gap decision theory
  integrated with energy management systems for isolated microgrids,''
  \emph{IEEE Trans. Ind. Inform.}, vol.~16, no.~8, pp. 5286--5299, Aug. 2020.

\bibitem{shuidian2024}
\BIBentryALTinterwordspacing
{China Renewable Energy Engineering Institute}, ``China renewable energy project
  cost management report 2024 edition,'' 2024. [Online]. Available:
  \url{https://h2.in-en.com/html/h2-2442597.shtml}
\BIBentrySTDinterwordspacing

{\color{black}\bibitem{davoodi2025scalable}
E.~Davoodi, F.~Capitanescu, M.~I. Alizadeh, and L.~Wehenkel, ``A scalable var
  planning methodology to mitigate reactive power scarcity during energy
  transition,'' \emph{IEEE Trans. Power Syst.}, pp. 1--11, 2025, Early Access.}

\bibitem{wu2023multi}
S.~Wu, J.~Lin, J.~Li, F.~Liu, Y.~Song, Y.~Xu, X.~Cheng, and Z.~Yu,
  ``Multi-timescale trading strategy for renewable power to ammonia virtual
  power plant in the electricity, hydrogen, and ammonia markets,'' \emph{IEEE
  Trans. Energy Mark. Policy Regul.}, vol.~1, no.~4, pp. 322--335, Dec. 2023.

{\color{black}\bibitem{huadian2024}
\BIBentryALTinterwordspacing
{China Huadian Corporation Ltd.}, ``Huadian {L}iaoning {T}ieling off-grid
  energy storage and hydrogen production integrated project,'' 2024. [Online].
  Available:
  \url{https://www.chd.com.cn/site/2/2024-05-08/966791b886fc4ba290fddff68d3422f4.html}
\BIBentrySTDinterwordspacing}

\bibitem{qiu2025dynamic}
Y.~Qiu, J.~Li, Y.~Zeng, Y.~Zhou, S.~Chen, X.~Qiu, B.~Zhou, G.~He, X.~Ji, and
  W.~Li, ``Dynamic operation and control of a multi-stack alkaline water
  electrolysis system with shared gas separators and lye circulation: A
  model-based study,'' \emph{arXiv preprint arXiv:2501.14576}, 2025.

\end{thebibliography}
\end{document}